\documentclass[a4paper,12pt]{article}

\usepackage{float}
\usepackage{amsmath} 
\usepackage{amssymb}
\usepackage{amsfonts}
\usepackage{epstopdf}
\usepackage{color,fancybox,graphicx}
\usepackage{url}
\usepackage{psfrag}
\usepackage{color}
\usepackage{algorithm, algorithmic}
\usepackage{pstricks}
\usepackage{graphicx,psfrag}
\usepackage{epsfig}

\usepackage{natbib}
\bibpunct{(}{)}{;}{a}{,}{,}

\numberwithin{equation}{section}

\newtheorem{lem}{Lemma}[section]
\newtheorem{cor}{Corollary}[section]
\newtheorem{pro}{Proposition}[section]
\newtheorem{theo}{Theorem}[section]

\newcommand{\bS}{\mathbf{S}}
\newcommand{\bs}{\mathbf{s}}
\newcommand{\bx}{\mathbf{x}}
\newcommand{\bX}{\mathbf{X}}

\newcommand{\bTheta}{\boldsymbol \Theta}
\newcommand{\btheta}{\boldsymbol \theta}
\newcommand{\bd}{\boldsymbol}

\begin{document}

\begin{center}
{\Large 
\textbf{\textsf{New Insights Into Approximate Bayesian Computation}}}
\medskip

\end{center}
{\bf G\'erard Biau\footnote{Corresponding author.}\\
{\it Universit\'e Pierre et Marie Curie\footnote{Research partially supported by the French National Research Agency under grant ANR-09-BLAN-0051-02 ``CLARA''.} \& Ecole Normale Sup{\'e}rieure\footnote{Research carried out within the INRIA project ``CLASSIC'' hosted by Ecole Normale Sup{\'e}rieure and CNRS.}, France}}\\
\textsf{gerard.biau@upmc.fr}
\bigskip

{\bf Fr\'ed\'eric C\'erou}\\
{\it INRIA Rennes -- Bretagne Atlantique, France}\\
\textsf{Frederic.Cerou@inria.fr}
\bigskip

{\bf Arnaud Guyader}\\
{\it Universit\'e Rennes 2 \& INRIA Rennes -- Bretagne Atlantique, France }\\
\textsf{arnaud.guyader@uhb.fr}
\medskip

\begin{abstract}
\noindent {\rm Approximate Bayesian Computation ({\sc abc} for short) is a family of computational techniques which offer an almost automated solution in situations where evaluation of the posterior likelihood is computationally prohibitive, or whenever suitable likelihoods are not available. In the present paper, we analyze the procedure from the point of view of $k$-nearest neighbor theory and explore the statistical properties of its outputs. We discuss in particular some asymptotic features of the genuine conditional density estimate associated with {\sc abc}, which is an interesting hybrid between a $k$-nearest neighbor and a kernel method.
\medskip

\medskip

}

\end{abstract}

\section{Introduction}
Let $\bd Y$ be a generic random observation which may, for example, take the form of a sample of independent and identically distributed (i.i.d.)~random variables. More generally, it may also be the first observations of a time series or a more complex random object, such as a {\sc dna} sequence. We denote by $\ell (\bd y | \btheta)$ the distribution (likelihood) of $\bd Y$, where $\btheta \in \mathbb R^p$ is an unknown parameter that we wish to estimate. In the Bayesian paradigm, the parameter itself is seen as a random variable $\bTheta $, and the likelihood $\ell(\bd y | \btheta)$ becomes the conditional distribution of $\bd Y$ given $\bTheta = \btheta$. The distribution $\pi(\btheta)$ of $\bTheta$ is called the prior distribution, while the distribution $\pi(\btheta|\bd y)$ of $\bTheta$ given $\bd Y=\bd y$ is termed posterior.
\medskip

When taking a Bayesian perspective, inference about the parameter $\bTheta$ typically proceeds via calculation or simulation of the posterior distribution $\pi(\btheta|\bd y)$. A variety of methods exist for inference in this context, such as rejection algorithms \citep[][]{Ripley}, Markov Chain
Monte Carlo ({\sc mcmc}) methods \citep[e.g., the Metropolis-Hastings algorithm,][]{Metropolis,Hastings}, and Importance Sampling \citep[][]{Ripley}. For a comprehensive introduction to the domain, the reader is referred to the monographs by \citet{Robert} and \citet{MR}. However, in some contexts, computation of the posterior is problematic, either because the size of the data makes the calculation computationally intractable, or because calculation
is impossible when using realistic models for how the data arises. Thus, despite their power and flexibility,  {\sc mcmc} procedures and their variants may prove irrelevant in a growing number of contemporary applications involving very large dimensions or complicated models. This computational burden typically arises in fields such as ecology, population genetics and image analysis, just to name a few.
\medskip

This difficulty has motivated a drive to more approximate approaches, in particular the
field of Approximate Bayesian Computation ({\sc abc} for short). In a nutshell,  {\sc abc} is a family of computational techniques which offer an almost automated solution in situations where evaluation of the likelihood is computationally prohibitive, or whenever suitable likelihoods are not available. The approach was originally mentioned, but not analyzed,
 by \citet{Rubin}. It was further developed in population genetics by \citet{Fu,Tavare,Pritchard,BZB}, who gave the name of Approximate Bayesian Computation to a family of likelihood-free inference methods. Since its original developments, the  {\sc abc} paradigm has successfully been applied to various scientific areas, ranging from archaeological science and ecology to epidemiology, stereology and protein network analysis. There are too many references to be included here, but the recent survey by \citet{Marin} offers  both a historical and  technical review of the domain. 
 \medskip
 
Before we go into more details on {\sc abc}, some more notation is required. We assume to be given a statistic $\bS$, taking values in $\mathbb R^m$. It is a function of the original observation $\bd Y$, with a dimension $m$ typically much smaller than the dimension of $\bd {Y}$. The statistic $\bS$ is supposed to admit a conditional density $f(\bs|\btheta)$ with respect to the Lebesgue measure on $\mathbb R^m$.  Note that, strictly speaking, we should write $\bS (\bd Y) $ instead of $\bS$. However, since there is no ambiguity, we continue to use the latter notation. As such, the statistic $\bS$ should be understood as a low-dimensional summary of $\bd Y$. It can be, for example, a sufficient statistic for the parameter $\bTheta $, but not necessarily. Assuming that $\bTheta$ is absolutely continuous with respect to the Lebesgue measure on $\mathbb R^p$, the conditional distribution of $\bTheta$ given $\bS = \bs$ has a density $g(\btheta | \bs) $ which, according to Bayes' rule, takes the form
$$g(\btheta|\bs)=\frac{f(\bs|\btheta)\pi(\btheta)}{\bar f(\bs)}, \quad \mbox{where }\bar f(\bs)=\int_{\mathbb R^p} f(\bs|\btheta)\pi(\btheta)\mbox{d}\btheta$$
is the marginal density of $\bS$. Finally, we denote by $\bd y_0$ the observed realization of $\bd Y$ (i.e., the data set), and let $\bs_0 (=\bs(\bd y_0))$ be the corresponding realization of $\bS$. Throughout the document, both $\bd y_0$ and $\bs_0$ should be considered as fixed quantities. 
\medskip

In its most common form, the generic {\sc abc} algorithm is framed as follows:
\begin{algorithm} 
\caption{Pseudo-code 1 of a generic {\sc abc} algorithm}
\label{algorithme1} 
\begin{algorithmic}
\REQUIRE A positive integer $N$ and a tolerance level $\varepsilon$.
\FOR{$i=1$ to $N$}
\STATE Generate $\btheta_i$ from the prior $\pi(\btheta)$;
\STATE Generate $\bd y_i$ from the likelihood $\ell(.|\btheta_i)$.
\ENDFOR
\RETURN The $\btheta_i$'s such that $\|\bs(\bd y_i)-\bs_0\|\leq \varepsilon$.
\end{algorithmic} 
\end{algorithm} 

The basic idea behind this formulation is that using a representative
enough summary statistic $\mathbf S$ coupled with a small enough
tolerance level $\varepsilon$ should produce a good approximation
of the posterior distribution. A moment's thought reveals that pseudo-code \ref{algorithme1} 
has the flavor of a nonparametric kernel conditional density estimation procedure, for which $\varepsilon$ 
plays the role of a bandwidth. This is, for example, the point of view that prevails in the analysis of \citet{Blum}, who explores the asymptotic bias and variance of kernel-type estimates of the posterior density $g(.|\bs_0)$ evaluated over the code outputs.
\medskip

However, as made transparent by \citet{Marin}, pseudo-code \ref{algorithme1}, despite its widespread diffusion, does not exactly match what people do in practice. A more accurate formulation is the following one:
\newpage
\begin{algorithm}[!h]
\caption{Pseudo-code 2 of a generic {\sc abc} algorithm}
\label{algorithme2} 
\begin{algorithmic}
\REQUIRE A positive integer $N$ and an integer $k_N$ between $1$ and $N$.
\FOR{$i=1$ to $N$}
\STATE Generate $\btheta_i$ from the prior $\pi(\btheta)$;
\STATE Generate $\bd y_i$ from the likelihood $\ell(.|\btheta_i)$.
\ENDFOR
\RETURN The $\btheta_i$'s such that $\bs(\bd y_i)$ is among the $k_N$-nearest neighbors of $\bs_0$.
\end{algorithmic} 
\end{algorithm} 
Algorithm \ref{algorithme1} and Algorithm \ref{algorithme2} are dual, in the sense that the number of accepted points is fixed in the second and random in the first, while their range is random in the second and fixed in the first. In practice, the parameter $N$ is chosen to be very large (typically of the order of $10^6$), while $k_N$ is most commonly expressed as a percentile. Thus, for example, the choice $N=10^6$ and a percentile $k_N/N=0.1\%$ allow to retain $1000$ simulated $\btheta_i$'s.
\medskip

From a nonparametric perspective, pseudo-code \ref{algorithme2} falls within the broad family of nearest neighbor-type procedures \citep[][]{Fix1,Loftsgaarden,Cover68}. Such procedures have the favor of practitioners, because they are fast, easy to compute and flexible. For implementation, they require only a measure of distance in the sample space, hence their popularity as a starting-point for refinement, improvement and adaptation to new settings \citep[see, e.g.,][Chapter 19]{DGL}. In any case, it is our belief that {\sc abc} should be analyzed in this context, and this is the point of view that is taken in the present article.
\medskip

\begin{figure}
\begin{center}
\input{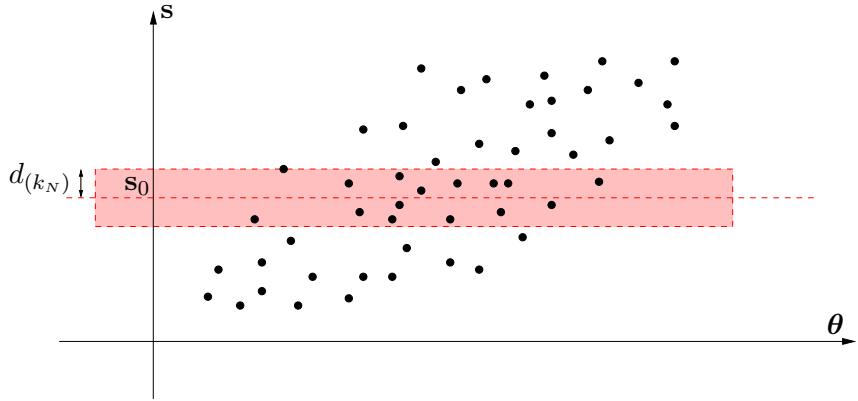}
\caption{Illustration of {\sc abc} in dimension $m=p=1$ ($d_{(k_N)}=\|\bS_{(k_N)}-\bs_0\|$).}
\label{figure1}
\end{center}
\end{figure}

In order to better understand the rationale behind Algorithm \ref{algorithme2}, denote by $(\bTheta_1, \bd Y_1), \hdots, (\bTheta_N,\bd Y_N)$ an i.i.d.~sample, with common joint distribution $\ell(\bd y|\btheta)\pi(\btheta)$. This sample is naturally associated with the i.i.d.~sequence $(\bTheta_1,\bS_1), \hdots, (\bTheta_N,\bS_N)$, where each pair has density $f(\bs|\btheta)\pi(\btheta)$. 
Finally, let $\bS_{(1)}, \hdots, \bS_{(k_N)}$ be the $k_N$-nearest neighbors of  $\bs_0$ among $\bS_1, \hdots, \bS_N$, and let $\bTheta_{(1)}, \hdots, \bTheta_{(k_N)}$ be the corresponding $\bTheta_i$'s (see Figure \ref{figure1} for an illustration in dimension $m=p=1$).

\medskip

With this notation, we see that the generic {\sc abc} Algorithm \ref{algorithme2} proceeds in two steps:
\begin{enumerate}
\item First, simulate (realizations of) an $N$-sample $(\bTheta_1,\bd Y_1), \hdots, (\bTheta_N,\bd Y_N)$;
\item Seconds, return (realizations of) the variables $\bTheta_{(1)}, \hdots, \bTheta_{(k_N)}$.
\end{enumerate}
This simple observation opens the way to a mathematical analysis of {\sc abc} via techniques based on nearest neighbors. In fact, despite a growing number of practical applications, theoretical results guaranteeing the validity of the approach are still lacking \citep[see][for results in this direction]{Wilkinson,Blum,Fearnhead}. Our present contribution is twofold:
\begin{enumerate}
\item[$(i)$] We offer in Section 2 an explicit result regarding the distribution of the algorithm outputs $(\bTheta_{(1)},\bS_{(1)}), \hdots, (\bTheta_{(k_N)},\bS_{(k_N)})$. Let $\mathcal B_m(\bs_0,\delta)$ denote the closed ball in $\mathbb R^m$ centered at $\bs_0$ with nonnegative radius $\delta$, i.e., $\mathcal B_m(\bs_0,\delta)=\{\bs \in \mathbb R^m:\|\bs-\bs_0\|\leq \delta\}$. In a nutshell, Proposition \ref{theorem1} reveals that, conditionally on the distance $d_{(k_N+1)}=\|\bS_{(k_N+1)}-\bs_0\|$, the simulated data set may be regarded as $k_N$ i.i.d.~realizations of the joint density of $(\bTheta,\bS)$ restricted to the cylinder $\mathbb R^p\times\mathcal B_m(\bs_0,d_{(k_N+1)})$. This result is important since it gives a precise description of the output distribution of {\sc abc} Algorithm \ref{algorithme2}.
\item[$(ii)$] For a fixed $\bs_0\in \mathbb R^m$, the estimate practitioners use most to infer the posterior density $g(.|\bs_0)$ at some point $\btheta_0 \in \mathbb R^p$ is
\begin{equation}
\label{Bruxelles}
\hat g_{N,\bs_0}(\btheta_0)=\frac{1}{k_Nh_N^p} \sum_{j=1}^{k_N}
K\left(\frac{\btheta_0-\bTheta_{(j)}}{h_N}\right),
\end{equation}
where $\{h_N\}$ is a sequence of positive real numbers (bandwidth) and $K$ is a nonnegative Borel measurable function (kernel) on $\mathbb R^p$. The idea is simple: In order to estimate the posterior, just look at the $k_N$-nearest neighbors of $\bs_0$ and smooth the corresponding $\bTheta_j$'s around $\btheta_0$. 
It should be noted that (\ref{Bruxelles}) is a smart hybrid between a $k$-nearest neighbor and a kernel density estimation procedure. It is different from the Rosenblatt-type \citep{Rosenblatt} kernel conditional density estimates proposed in \citet{BZB} and further explored by \citet{Blum}. In Section 3 and Section 4, we establish some consistency properties of this genuine estimate and discuss its rates of convergence. 
\end{enumerate}
For the sake of clarity, proofs are postponed to Section 5 and Section 6. An appendix at the end of the paper offers some new results on convolution and approximation of the identity.
\medskip

To conclude this introduction, we would like to make a few comments on the topics that will {\bf not} be addressed in the present document. An important part of the performance of the {\sc abc} approach, especially for high-dimensional data sets, relies upon a good choice of the summary statistic $\bS$. In many practical applications, this statistic is picked by an expert in the field, without any particular guarantee of success. A systematic approach to choosing such a statistic, based upon a sound theoretical
framework, is currently under active investigation in the Bayesian community. This important issue will not be pursued further here. As a good starting point, the interested reader is referred to \citet{Joyce}, who develop a sequential scheme for scoring statistics according to whether their inclusion in the analysis will substantially improve the quality of inference. Similarly, we will not address issues regarding how to enhance efficiency of {\sc abc} and its variants, as for example with the sequential techniques of \citet{Sisson} and \citet{BCMR}. Nor won't we explore the important question of {\sc abc} model choice, for which theoretical arguments are still missing \citep[][]{RCMP,MPRR}.
\section{Distribution of {\sc abc} outputs}
We continue to use the notation of Section 1 and recall in particular that $(\bTheta_1,\bS_1), \hdots, (\bTheta_N,\bS_N)$ are i.i.d.~$\mathbb R^p\times \mathbb R^m$-valued random variables, with common probability density $f(\btheta,\bs)=f(\bs|\btheta)\pi(\btheta)$. Both $\mathbb R^p$ (the space of $\bTheta_i$'s) and $\mathbb R^m$ (the space of $\bS_i$'s) are equipped with the Euclidean norm $\|.\|$. In this section, attention is focused on analyzing the distribution of the algorithm outputs $(\bTheta_{(1)},\bS_{(1)}), \hdots, (\bTheta_{(k_N)},\bS_{(k_N)})$.
\medskip

In what follows, we keep $\bs_0$ fixed and denote by $d_{i}$ the (random) distance between $\bs_0$ and $\bS_i$. (To be rigorous, we should write $d_{i}(\bs_0)$, but since no confusion can arise we write it simply $d_{i}$.) Similarly, we let $d_{(i)}$ be the distance between $\bs_0$ and its $i$-th nearest neighbor among $\bS_1, \hdots, \bS_N$, that is
$$d_{(i)}=\|\bS_{(i)}-\bs_0\|.$$
(If distance ties occur, a tie-breaking strategy must be defined. For example, if $\|\bS_i-\bs_0\|=\|\bS_j-\bs_0\|$, $\bS_i$ may be declared ``closer'' if $i <j$, i.e., the tie-breaking is done by indices. Note however that ties occur with probability 0 since all random variables are absolutely continuous.) It is assumed throughout the paper that $N\geq 2$ and $1\leq k_N\leq N-1$. 
\medskip

Rearranging the $k_N$ (ordered) statistics $(\bTheta_{(1)},\bS_{(1)}), \hdots, (\bTheta_{(k_N)},\bS_{(k_N)})$ in the original order of their outcome, one obtains the $k_N$ (non-ordered) random variables $(\bTheta^{\star}_{1},\bS^{\star}_{1}), \hdots, (\bTheta^{\star}_{k_N},\bS^{\star}_{k_N})$.  Our first result is concerned with the conditional distributions
$$\mathcal L\left \{ (\bTheta^{\star}_1,\bS^{\star}_1), \hdots,(\bTheta^{\star}_{k_N},\bS^{\star}_{k_N})  \,|\, d_{(k_N+1)}\right\}$$
and
$$\mathcal L \left \{ (\bTheta_{(1)},\bS_{(1)}), \hdots,(\bTheta_{(k_N)},\bS_{(k_N)})\,|\, d_{(k_N+1)}\right\}.$$
Recall that the collection of all $\mathbf s_0 \in \mathbb R^m$ with $\int_{\mathcal B_m(\mathbf s_0,\delta)}\bar f(\bs)\mbox{d}\bs>0$ for all $\delta>0$ is called the support of $\bar f$. 
\begin{pro}[Distribution of {\sc abc} outputs]
\label{theorem1}
Assume that $\bs_0$ belongs to the support of $\bar f$. Let $(\tilde\bTheta_1,\tilde \bS_1), \hdots,(\tilde\bTheta_{k_N},\tilde
\bS_{k_N})$ be i.i.d.~random variables, with common probability density (conditional on $d_{(k_N+1)}$) 
\begin{equation}
\label{ladensite}
\frac{{\mathbf 1}_{[ \|
  \bs-\bs_0\|\leq d_{(k_{N}+1)}]} f(\btheta,\bs)}{\displaystyle \int_{\mathbb R^p}\int_{\mathcal B_m(\bs_0,d_{(k_N+1)})}f(\btheta,\bs)\emph{d}{\btheta}\emph{d}{\bs}}.
  \end{equation}
 Then
$$ \mathcal L \left \{ (\bTheta^{\star}_1,\bS^{\star}_1), \hdots,(\bTheta^{\star}_{k_N},\bS^{\star}_{k_N}) \,|\, d_{(k_N+1)}\right\}=\mathcal L\left \{(\tilde\bTheta_1,\tilde \bS_1), \hdots,(\tilde\bTheta_{k_N},\tilde
\bS_{k_N})\right\}.
$$
Moreover
\begin{align*}
&\mathcal L \left \{ (\bTheta_{(1)},\bS_{(1)}), \hdots,(\bTheta_{(k_N)},\bS_{(k_N)})\,|\, d_{(k_N+1)}\right\}\\
&\quad =\mathcal L\left \{(\tilde\bTheta_{(1)},\tilde \bS_{(1)}), \hdots,(\tilde\bTheta_{(k_N)},\tilde
\bS_{(k_N)})\right\}.
\end{align*}
\end{pro}
Note, since $\bs_0$ belongs by assumption to the support of $\bar f$, that the normalizing constant in the denominator of (\ref{ladensite}) is positive. This theorem may be regarded as an extension of a result of \citet{Kaufmann}, who provide explicit representations of the conditional distribution of an empirical point process given some order statistics. However, the present Bayesian setting is not covered by the conclusions of \citet{Kaufmann}, and our proof actually relies on simpler arguments.
\medskip

The main message of Proposition \ref{theorem1} is that, {\bf conditionally on $d_{(k_N+1)}$}, one can consider the $k_N$-tuple $(\bTheta_{(1)},\bS_{(1)}), \hdots,(\bTheta_{(k_N)},\bS_{(k_N)})$ as an ordered sample drawn according to the probability density (\ref{ladensite}). Alternatively, the (unordered) simulated values may be treated like i.i.d.~realizations of  variables with common density proportional to ${\mathbf 1}_{[ \| \bs-\bs_0\|\leq d_{(k_{N}+1)}]} f(\btheta,\bs)$. Conditionally on $d_{(k_N+1)}$, the accepted $\btheta_j$'s are nothing but i.i.d.~realizations of the probability density
$$\frac{\displaystyle \int_{\mathcal B_m(\bs_0,d_{(k_N+1)})}f(\btheta,\bs)\mbox{d}{\bs}}{\displaystyle \int_{\mathbb R^p}\int_{\mathcal B_m(\bs_0,d_{(k_N+1)})}f(\btheta,\bs)\mbox{d}{\btheta}\mbox{d}{\bs}}.$$
Although this conclusion is intuitively clear, its proof requires a careful mathematical analysis.
\medskip

As will be made transparent in the next section, Proposition \ref{theorem1} plays a key role 
in the mathematical analysis of the natural conditional density estimate associated with {\sc abc}
methodology. In fact, investigating {\sc abc} in terms of nearest neighbors has other important consequences. Suppose, for example, that we are interested in estimating some finite conditional expectation
$\mathbb E[\varphi(\bTheta)|\bS=\bs_0]$, where the random variable $\varphi(\bTheta)$ is bounded. This includes in particular the important setting where $\varphi$ is polynomial and one wishes to estimate the conditional moments of $\bTheta$. Then, provided $k_N/\log\log N \to \infty$ and $k_N/N
\to 0$ as $N\to \infty$, it can be shown that for almost all $\bs_0$ (with respect to the distribution of $\bS$), with
probability 1, 
\begin{equation}
\label{gold}
\frac{1}{k_N}\sum_{j=1}^{k_N} \varphi\left (\bTheta_{(j)}\right) \to 
\mathbb E[\varphi(\bTheta)|\bS=\bs_0].
\end{equation}
Proof of such a result uses the full power of the vast and rich nearest neighbor estimation theory. To be more precise, let us make a quick detour through this theory and consider an i.i.d. sample $(\bX_1,Z_1),\dots,(\bX_N,Z_N)$ taking values in $\mathbb R^m \times \mathbb
R$, where the output variables $Z_i$'s are bounded. Assume, to keep things simple, that the $\bX_i$'s have a probability density and that our goal is to assess the regression function $r(\bx) = \mathbb
E[Z\,|\,\bX=\bx]$, $\bx \in \mathbb R^m$. In this context, the $k$-nearest neighbor regression function estimate of $r$ \citep[][]{Royall,Cover68,Stone77} takes the form
$$\hat r_N(\bx) = \frac{1}{k_N} \sum_{j=1}^{k_N} Z_{(j)}, \quad \bx \in \mathbb R^m,$$ 
where $Z_{(j)}$ is the $Z$-observation corresponding to $\bX_{(j)}$, the $j$-th-closest point to $\bx$ among $\bX_1, \hdots, \bX_N$.  Denoting by $\mu$ the distribution of $\bX_1$, it is proved in Theorem 3 of \citet{devroye82} that provided $k_N/\log\log N \to \infty$ and $k_N/N \to 0$, for $\mu$-almost all $\bx$,
$$\hat r_N(\bx)\to r(\bx)\quad \mbox{with probability 1 as } N\to \infty.$$ 
This result can be transposed without further effort to our {\sc abc} setting via the
correspondence $\varphi(\bTheta) \leftrightarrow Z$ and $\bS \leftrightarrow \bX$, thereby establishing validity of (\ref{gold}). The decisive step towards that conclusion is accomplished by making a connection between {\sc abc} and nearest neighbor methodology. We leave it to the reader to draw his own conclusions as to further possible utilizations of this correspondence.
 
\section{Mean square error consistency}

As in Section 2, we keep the conditioning vector $\bs_0$ fixed and consider the i.i.d.~sample $(\bTheta_1,\bS_1), \hdots, (\bTheta_N,\bS_N)$, where each pair is distributed according to the probability density $f(\btheta,\bs)=f(\bs|\btheta)\pi(\btheta)$ on $\mathbb R^p\times \mathbb R^m$.  Based on this sample, our new objective is to estimate the posterior density $g(\btheta_0|\bs_0)$, $\btheta_0\in \mathbb R^p$. 
This estimation step is an important ingredient of the Bayesian analysis, whether this may be for visualization purposes or more involved mathematical achievements.
\medskip

As exposed in the introduction, the natural {\sc abc}-companion estimate of $g(\btheta_0|\bs_0)$ takes the form
\begin{equation}
\label{Nougaro}
\hat g_N(\btheta_0)=\frac{1}{k_Nh_N^p} \sum_{j=1}^{k_N}
K\left(\frac{\btheta_0-\bTheta_{(j)}}{h_N}\right),\quad \btheta_0 \in \mathbb R^p,
\end{equation}
where $\{h_N\}$ is a sequence of positive real numbers (bandwidth) and $K$ is a nonnegative Borel measurable function (kernel) on $\mathbb R^p$. (To reduce the notational burden, we dropped the dependency of the estimate upon $\bs_0$, keeping in mind that $\bs_0$ is held fixed.) Kernel estimates were originally studied in density estimation by \citet{Rosenblatt} and \citet{Parzen}, and were latter introduced in regression estimation by \citet{NadarayaReg1, NadarayaReg2} and \citet{Watson64}. The origins of $k$-nearest neighbor density estimation go back to \citet{Fix1} and \citet{Loftsgaarden}. Kernel estimates have been extended to the conditional density setting by \citet{Rosenblatt}, who proceeds by separately inferring the bivariate density $f(\btheta,\bs)$ of $(\bTheta,\bS)$ and the marginal density of $\bS$. Rosenblatt's estimate reads
$$\tilde g_N(\btheta_0)=\frac{\sum_{i=1}^NL\left(\frac{\bs_0-\bS_{i}}{\delta_N}\right)K\left(\frac{\btheta_0-\bTheta_{i}}{h_N}\right)}{h_N^p\sum_{i=1}^NL\left(\frac{\bs_0-\bS_{i}}{\delta_N}\right)},$$
where $L$ is a kernel in $\mathbb R^m$, and $\delta_N$ is the corresponding bandwidth. {\sc abc}-compatible estimates of this type have been discussed in \citet{BZB} and further explored by \citet{Blum} (\citealp[additional references for the conditional density estimation problem are][]{Hyndman,GK,Faugeras}, and the survey of \citealp[][]{Hansen}).
\medskip

The conditional density estimate we are interested in is different, in the sense that it has both the flavor of a $k$-nearest neighbor approach (it retains only the $k_N$-nearest neighbors of $\bs_0$ among $\bS_1, \hdots, \bS_N$) and a kernel method (it smoothes the corresponding $\bTheta_j$'s). Obviously, the main advantage of (\ref{Nougaro}) over its kernel-type competitors is its simplicity (it does not involve evaluation of a ratio, with a denominator that can be small), which makes it easy to implement. 
\medskip

A related procedure to density estimation has been originally proposed by \citet{brmepu77}, who suggested varying the kernel bandwidth with respect to the sample points. Various extensions and modifications of the  \citet{brmepu77} estimate have been later proposed in the literature.  The rationale behind the approach is to combine the desirable smoothness properties of kernel estimates with the data-adaptive character of nearest neighbor procedures. Particularly influential papers in the study of variable kernel estimates were those of \citet{Abramson} and \citet{HM}, who showed how variable bandwidths with positive kernels can nevertheless induce convergence rates usually attainable with fixed bandwidths and fourth order kernels. For a complete and comprehensive description of variable kernel estimates and their properties, we refer the reader to \citet{Jones}.
\medskip

Our goal in this section is to investigate some consistency properties of the {\sc abc}-companion estimate (\ref{Nougaro}). Pointwise mean square error consistency is proved in Theorem \ref{theorem2} and mean integrated square error consistency is established in Theorem \ref{theorem3}. We stress that this part of the document is concerned with minimal conditions of convergence. We did indeed try to reduce as much as possible the assumptions on the various unknown probability densities by resorting to real analysis arguments.
\medskip

The following assumptions on the kernel will be needed throughout the paper:
\medskip

\noindent\textbf{Assumption [K1]}\quad The kernel $K$ is nonnegative and belongs to $L^1(\mathbb R^p)$, with 
$\int_{\mathbb  R^p} K(\btheta)\mbox{d}\btheta=1$. Moreover, the function $\sup_{\|\mathbf y\|\geq \|\btheta\|} |K(\mathbf y)|$, $\btheta \in \mathbb R^{p}$, is in $L^1(\mathbb R^{p})$.
\medskip

Assumption set $[{\bf K1}]$ is in no way restrictive and is satisfied by all standard kernels such as, for example, the naive kernel
$$K(\btheta)=\frac{1}{V_p}\mathbf 1_{\mathcal B_p(\mathbf 0,1)}(\btheta),$$
where $V_p$ is the volume of the closed unit ball $\mathcal B_p(\mathbf 0,1)$ in $\mathbb R^p$,
or the Gaussian kernel
$$K(\btheta)=\frac{1}{(2\pi)^{p/2}} \exp\left (-\|\btheta\|^2/2\right).$$
We recall for further references that, in the $p$-dimensional Euclidean space,
\begin{equation}
\label{volume}
V_p=\frac{\pi^{p/2}}{\Gamma \left(1+\frac{p}{2}\right)},
\end{equation}
where $\Gamma(.)$ is the gamma function. Everywhere in the document, we denote by $\lambda_p$ (respectively, $\lambda_m$) the Lebesgue measure on $\mathbb R^p$ (respectively, $\mathbb R^m$) and set, for any positive $h$,
$$K_{h}(\btheta)=\frac{1}{h^p}K(\btheta/h), \quad \btheta \in \mathbb R^p.$$
We note once and for all that, under Assumption $[{\bf K1}]$, $\int_{\mathbb R^p} K_h(\btheta)\mbox{d}\btheta=1$.
\medskip

The first crucial result from real analysis that is needed here is the so-called Lebesgue's differentiation theorem \citep[see, e.g., Theorem 7.16 in][]{Wheeden}, which asserts that if $\varphi$ is a locally integrable function in $\mathbb R^n$, then
$$\frac{1}{V_n\delta^n} \int_{\mathcal B_n(\bx_0,\delta)}\left | \varphi(\bx)-\varphi(\bx_0)\right|\mbox{d}\bx\to 0\quad \mbox{as } \delta\to 0$$
for $\lambda_n$-almost all $\bx_0 \in \mathbb R^n$. A point $\bx_0$ at which this statement is valid is called a Lebesgue point of $\varphi$. In the proofs, we shall in fact need some convolution-type variations around the Lebesgue's theorem regarding the prior density $\pi$. These important results are gathered in the next theorem, whose proof can be found in Theorem 1, page 5 and Theorem 2, pages 62-63 of \citet{stein}.
\begin{theo}
\label{WZ1}
Let $K$ be a kernel satisfying Assumption $[{\bf K1}]$, and let the function $\pi^{\star}$ be defined on $\mathbb R^p$ by
$$\btheta_0 \mapsto {\pi}^{\star}(\btheta_0)=\sup_{h>0} \left [\int_{\mathbb R^p}K_{h}(\btheta_0-\btheta) \pi(\btheta)\emph{d}\btheta\right].$$
\begin{enumerate}
\item[$(i)$] For $\lambda_p$-almost all $\btheta_0 \in \mathbb R^p$, one has
$$ \int_{\mathbb R^p} K_h(\btheta_0-\btheta)\pi(\btheta)\emph{d}\btheta\to\pi(\btheta_0)\quad \mbox{as } h\to 0.$$
\item[$(ii)$] The quantity ${\pi}^{\star}(\btheta_0)$ is finite for $\lambda_p$-almost all $\btheta_0 \in \mathbb R^p$.
\item[$(iii)$] For any $q>1$, the function ${\pi}^{\star}$ is in $L^q(\mathbb R^p)$ whenever $\pi$ is in $L^q(\mathbb R^p)$.
\end{enumerate}
\end{theo}
When $K$ is chosen to be the naive kernel, the function $\pi^{\star}$ of Theorem \ref{WZ1} is called the Hardy-Littlewood maximal function of $\pi$. It should be understood as a gauge of the size of the averages of $\pi$ around $\btheta_0$. 
\medskip

We shall also need an equivalent of Theorem \ref{WZ1} for the joint
density $f$, which this time is defined on $\mathbb R^p \times \mathbb
R^m$. Things turn out to be slightly more complicated in this case if
one is willing pairs of points $(\btheta_0,\bs_0)$ to be approached as
$(h,\delta)\to (0,0)$ by general product kernels over $\mathbb
R^p\times \mathbb R^m$. These kernels take the form $K_h(.)\otimes
L_{\delta}(.)$, without any restriction on the joint behavior of $h$
and $\delta$ (in particular, we do not impose that $h=\delta$). The
so-called Jessen-Marcinkiewicz-Zygmund theorem (\citealp[][]{JMZ}, see also
\citealp[][Chapter 17, pages 305-309]{Zygmund}) answers the
question for naive kernels, at the price of a slight integrability
assumption on $f$. On the other hand, the literature offers
surprisingly little help for general kernels, with the exception of arguments presented in \citet{DK}. This is
astonishing since this real analysis issue is at the basis of
pointwise convergence properties of multivariate kernel estimates and
indeed most density estimates. To fill the gap, we begin with the
following theorem, which is tailored to our {\sc abc} context (that is, when the second kernel $L$ is restricted to be the naive one). A more general result (that is, for both $K$ and $L$ general kernels) together with interesting new results on convolution and approximation of the identity are given in the Appendix section, at the end
of the paper (Theorem \ref{WZ2} is thus a consequence of Theorem \ref{proannex2}). In the sequel, notation $u^+$ means $\max(u,0)$.
\begin{theo}
\label{WZ2}
Let $K$ be a kernel satisfying Assumption $[{\bf K1}]$, and let the function $f^{\star}$ be defined on $\mathbb R^p\times \mathbb R^m$ by
$$(\btheta_0,\bs_0)\mapsto f^{\star}(\btheta_0,\bs_0)=\sup_{h>0,\delta>0} \left [\frac{1}{V_m\delta^m}\int_{\mathbb R^p}\int_{\mathcal B_m(\mathbf s_0,\delta)}K_{h}(\btheta_0-\btheta)  f(\btheta,\bs)\emph{d}\btheta\emph{d}\bs\right].$$
\begin{enumerate}
\item[$(i)$] If 
\begin{equation}
\label{NellyFurtado}
\int_{\mathbb R^p}\int_{\mathbb R^m} f(\btheta,\bs) \log^{+}f(\btheta,\bs)\emph{d}\btheta\emph{d}\bs<\infty
\end{equation}
then, for $\lambda_p \otimes \lambda_m$-almost all $(\btheta_0,\bs_0)\in \mathbb R^p\times \mathbb R^m$,
$$\lim_{(h,\delta) \to (0,0)}\frac{1}{V_m\delta^m}\int_{\mathbb R^p}\int_{\mathcal B_m(\mathbf s_0,\delta)}K_{h}(\btheta_0-\btheta)  f(\btheta,\bs)\emph{d}\btheta\emph{d}\bs= f(\btheta_0,\bs_0).$$
\item[$(ii)$] If condition (\ref{NellyFurtado}) is satisfied, then $f^{\star}(\btheta_0,\bs_0)$ is finite for $\lambda_p \otimes \lambda_m$-almost all $(\btheta_0,\bs_0)\in \mathbb R^p\times \mathbb R^m$.
\item[$(iii)$] For any $q>1$, the function $f^{\star}$ is in $L^q(\mathbb R^p\times \mathbb R^m)$ whenever $f$ is in $L^q(\mathbb R^p\times \mathbb R^m)$.
\end{enumerate}
\end{theo}
A remarkable feature of Theorem \ref{WZ2} $(i)$ is that the result is true as soon as $(h,\delta) \to (0,0)$, without any restriction on these parameters. This comes however at the price of the mild integrability assumption (\ref{NellyFurtado}), which is true, in particular, if $f$ is in any $L^q(\mathbb R^p\times \mathbb R^m)$, $q>1$.
\medskip

Recall that we denote by $\bar f$ the marginal density of $f(\btheta,\bs)$ in $\bs$, that is
$$\bar f(\bs)=\int_{\mathbb R^p} f(\btheta,\bs)\mbox{d}\btheta, \quad \bs \in \mathbb R^m.$$
We are now in a position to state the two main results of this section.
 \begin{theo}[Pointwise mean square error consistency]
\label{theorem2}
Assu\-me th\-at the kernel $K$ is bounded and satisfies Assumption $[{\bf K1}]$. Assume, in addition, that the joi\-nt probability density $f$ is such that
$$\int_{\mathbb R^p}\int_{\mathbb R^m} f(\btheta,\bs) \log^{+} f(\btheta,\bs)\emph{d}\btheta\emph{d}\bs<\infty.$$
Then, for $\lambda_p \otimes \lambda_m$-almost all $(\btheta_0,\bs_0)\in \mathbb R^p\times \mathbb R^m$, with $\bar f(\bs_0)>0$, if $k_N\to \infty$, $k_N/N\to 0$, $h_N\to 0$ and $k_Nh_N^p\to \infty$,
$$\mathbb E \left[ \hat g_N(\btheta_0) -g(\btheta_0|\bs_0)\right]^2 \to 0 \quad \mbox{as } N\to \infty.$$
\end{theo}
It is stressed that the integral assumption required on $f$ is mild. It is for example satisfied whenever $f$ is bounded from above or whenever $f$ belongs to $L^q(\mathbb R^p\times \mathbb R^m)$ with $q>1$. There are, however, situations where this assumption is not satisfied. As an illustration, take $p=m=1$ and let 
$${\cal T}=\left\{(\btheta,\bs)\in\mathbb R\times \mathbb R:\btheta>0, \bs>0, \btheta+\bs\leq\frac{1}{2}\right\}.$$
Clearly,
$$\iint_{{\cal T}}\frac{1}{(\btheta+\bs)^2\log^2(\btheta+\bs)}\mbox{d}\btheta\mbox{d}\bs<\infty.$$
Choose
$$f(\btheta,\bs)=\frac{C}{(\btheta+\bs)^2\log^2(\btheta+\bs)}\mathbf 1_{[(\btheta,\bs)\in{\cal T}]},$$
where $C$ is a normalizing constant ensuring that $f$ is a probability density.
Then
$$\int_{\mathbb R^p}\int_{\mathbb R^m} f(\btheta,\bs) \mbox{d}\btheta\mbox{d}\bs=1$$
whereas
$$\int_{\mathbb R^p}\int_{\mathbb R^m} f(\btheta,\bs) \log^{+} f(\btheta,\bs)\mbox{d}\btheta\mbox{d}\bs=\infty.$$
\medskip

Theorem \ref{theorem3} below states that the estimate $\hat g_N$ is also consistent with respect to the mean integrated square error criterion.
\begin{theo}[Mean integrated square error consistency]
\label{theorem3}
Assume th\-at the kernel $K$ belongs to $L^2(\mathbb R^p)$ and satisfies Assumption $[{\bf K1}]$. Assume, in addition, that the joint probability density $f$ and the prior $\pi$ are in $L^2(\mathbb R^p\times \mathbb R^m)$ and $L^2(\mathbb R^p)$, respectively. Then, for $\lambda_m$-almost all $\bs_0\in \mathbb R^m$, with $\bar f(\bs_0)>0$, if $k_N\to \infty$, $k_N/N\to 0$, $h_N\to 0$ and $k_Nh_N^p\to \infty$,
$$\mathbb E\left [ \int_{\mathbb R^p}\left[ \hat g_N(\btheta_0) -g(\btheta_0|\bs_0)\right]^2\emph{d}\btheta_0\right] \to 0 \quad \mbox{as } N\to \infty.$$
\end{theo}
Here again, the regularity assumptions required on $f$ and $\pi$ are minimal. One could envisage an additional degree of smoothing in the estimate (\ref{Nougaro}) by observing that taking the $k_N$ nearest neighbors of $\bs_0$ can be viewed as the uniform kernel case of the more general quantity 
$$\sum_{i=1}^N L \left (\frac{\bs_0-\bS_i}{\|\bS_{(k_N)}-\bs_0\|}\right),$$
which allows unequal weights to be given to the $\bS_i$'s. The corresponding smoothed conditional density estimate is defined by
$$\tilde g_{N}(\btheta_0)=\frac{\sum_{i=1}^{N}L \left (\frac{\bs_0-\bS_i}{\|\bS_{(k_N)}-\bs_0\|}\right)K\left(\frac{\btheta_0-\bTheta_{i}}{h_N}\right)}{h_N^p \sum_{i=1}^N L \left (\frac{\bs_0-\bS_i}{\|\bS_{(k_N)}-\bs_0\|}\right)}.$$
Thus, $\hat g_N$ is the uniform kernel case of $\tilde g_N$.  The asymptotic properties of $\tilde g_N$, which are beyond the scope of the present article, will be explored elsewhere by the authors. A good starting point are the papers by \citet{MY,Moore2} and \citet{MackRosen}, who study various properties of similar kernel-type nearest neighbor procedures for density estimation.

\section{Rates of convergence}

In this section, we go one step further in the analysis of the {\sc abc}-companion estimate $\hat g_N$ by studying its mean integrated square error rates of convergence. We follow the notation of Section 3 and try to keep the assumptions on unknown mathematical objects as mild as possible. Introduce the multi-index notation
$$|\beta|=\beta_1+\hdots+\beta_n, \quad \beta!=\beta_1!\hdots\beta_n!,\quad \bx^{\beta}=x_1^{\beta_1}\hdots x_n^{\beta_n}$$
for $\beta=(\beta_1, \hdots, \beta_n) \in \mathbb N^{n}$ and $\bx \in \mathbb R^n$. If all the $k$-order derivatives of some function $\varphi:\mathbb R^n\to \mathbb R$ are continuous at $\mathbf x_0\in \mathbb R^n$ then, by Schwarz's theorem, one can change the order of mixed derivatives at $\mathbf x_0$, so the notation
$$D^{\beta}\varphi(\mathbf x_0)=\frac{\partial^{|\beta|}\varphi(\mathbf x_0)}{\partial x_1^{\beta_1}\hdots\partial x_n^{\beta_n}}, \quad |\beta|\leq k$$
for the higher-order partial derivatives is justified in this situation.
\medskip

In the sequel, we shall need the following sets of
assumptions. Recall that the collection of all $\mathbf s_0 \in \mathbb R^m$ with $\int_{\mathcal B_m(\mathbf s_0,\delta)}\bar f(\bs)\mbox{d}\bs>0$ for all $\delta>0$ is called the support of $\bar f$. 
\medskip

\noindent\textbf{Assumption [A1]}\quad The marginal probability density $\bar f$ has compact support with diameter $L>0$ and is three times continuously differentiable.
\medskip

\noindent\textbf{Assumption [A2]}\quad The joint probability density $f$ is in $L^2(\mathbb R^p\times \mathbb R^m)$. Moreover, for fixed $\bs_0$, the functions
\begin{align*}
\btheta_0 &\mapsto \frac{\partial^2 f(\btheta_0,\bs_0)}{\partial \theta_{i_1}\partial \theta_{i_2}}, \quad 1\leq i_1,i_2\leq p\\
\mbox{and }\btheta_0 & \mapsto \frac{\partial^2 f(\btheta_0,\bs_0)}{\partial s_{j}^2}, \quad 1\leq j\leq m
\end{align*}
are defined and belong to $L^2(\mathbb R^p)$.
\medskip

\noindent\textbf{Assumption [A3]}\quad The joint probability density $f$  is three times continuously differentiable on $\mathbb R^p \times \mathbb R^m$ and, for any multi-index $\beta$ satisfying $|\beta|=3$,
$$\sup_{\bs\in\mathbb R^m} \int_{\mathbb R^p}\left[D^\beta f(\btheta,\bs)\right]^2\mbox{d}\btheta<\infty.$$ 
\medskip

It is also necessary to put some mild additional restrictions on the kernel.
\medskip

\noindent\textbf{Assumption [K2]}\quad The kernel $K$ is symmetric and belongs to $L^2(\mathbb R^p)$. Moreover,  for any multi-index $\beta$ satisfying $|\beta|\in\{1,2,3\}$,
$$\int_{\mathbb R^p} \left|\btheta^\beta\right| K(\btheta)\mbox{d}\btheta<\infty.$$
\medskip

We finally define
$$\xi_0=\inf_{0<\delta\leq L}\frac{1}{\delta^m} \int_{\mathcal B_m(\bs_0,\delta)}\bar f(\bs)\mbox{d}{\bs},$$
and introduce the following quantities, which are related to the average distance between $\bs_0$ and its $k_N$-th nearest neighbor (see Proposition \ref{distances} and Proposition \ref{distancesbis}):
\begin{align*}
D_m(k_N)&=\frac{m}{\xi_0^{2/m}(m-2)}\left (\frac{k_N+1}{N+1}\right)^{2/m} -\frac{L^{2-m}}{\xi_0(m/2-1)}\frac{k_N+1}{N+1},\\
\Delta_m(k_N)&=\frac{m}{\xi_0^{4/m}(m-4)}\left (\frac{k_N+1}{N+1}\right)^{4/m} -\frac{L^{4-m}}{\xi_0(m/4-1)}\frac{k_N+1}{N+1},\\
D(k_N)&=\frac{1}{\xi_0}\left(1+\log \left (\xi_0 L^2 \frac{N+1}{k_N+1}\right)\right)\frac{k_N+1}{N+1},\\
\Delta(k_N)&=\frac{1}{\xi_0}\left(1+\log \left (\xi_0 L^4 \frac{N+1}{k_N+1}\right)\right)\frac{k_N+1}{N+1}.
\end{align*}

The next theorem makes precise the mean integrated square error rates of convergence of $\hat g_N(.)$ towards $g(.|\bs_0)$. 
\begin{theo}
\label{theorem4}
Let $K$ be a kernel satisfying assumptions $[{\bf K1}]$ and  $[{\bf K2}]$. Let $\bs_0$ be a Lebesgue point of $\bar f$ such that $\bar f(\bs_0)>0$. Assume that Assumptions $[{\bf A1}]$-$[{\bf A3}]$ are satisfied. Then, letting
\begin{align*}
\phi_1(\btheta_0,\bs_0)&=\frac{1}{2}\sum_{i_1,i_2=1}^p\frac{\partial^2 f(\btheta_0,\bs_0)}{\partial \theta_{i_1}\partial \theta_{i_2}}\int_{\mathbb R^p} \theta_{i_1}\theta_{i_2}K (\btheta)\emph{d}\btheta\\
\phi_2(\btheta_0,\bs_0)&=\frac{1}{2m+4}\sum_{j=1}^m\frac{\partial^2 f(\btheta_0,\bs_0)}{\partial s_j^2},\\
\phi_3(\bs_0)&=\frac{1}{2m+4}\sum_{j=1}^m\frac{\partial^2 \bar f(\bs_0)}{\partial s_j^2},
\end{align*}
and
\begin{align*}
\Phi_1(\bs_0)& =\frac{1}{{\bar f}^2(\bs_0)}\int_{\mathbb R^p} \phi_1^2(\btheta_0,\bs_0)\emph{d}\btheta_0,\\
\Phi_2(\bs_0)&=\frac{1}{{\bar f}^4(\bs_0)}\int_{\mathbb R^p} \left [\phi_2(\btheta_0,\bs_0)\bar f(\bs_0)-\phi_3(\bs_0)f(\btheta_0,\bs_0)\right]^2\emph{d}\btheta_0,\\
\Phi_3(\bs_0) &=\frac{2}{{\bar f}^3(\bs_0)}\int_{\mathbb R^p}  \phi_1(\btheta_0,\bs_0)\left [\phi_2(\btheta_0,\bs_0)\bar f(\bs_0)-\phi_3(\bs_0)f(\btheta_0,\bs_0)\right]\emph{d}\btheta_0,
\end{align*}
one has:
\begin{enumerate}
\item {\bf For $\boldsymbol{m=2}$},
\begin{align*}
&\mathbb E\left [ \int_{\mathbb R^p}\left[ \hat g_N(\btheta_0) -g(\btheta_0|\bs_0)\right]^2\emph{d}\btheta_0\right]\\
& \quad =\left(\Phi_1(\bs_0)h_N^4+\Phi_2(\bs_0)\Delta_2(k_N)+\Phi_3(\bs_0)h_N^2D(k_N)+\frac{\int_{\mathbb R^p} K^2(\btheta)\emph{d}\btheta}{k_Nh_N^{p}}\right)\\
& \qquad \times\left (1+\emph{o}(1)\right).
\end{align*}
\item {\bf For $\boldsymbol{m=4}$},
$$$$
\begin{align*}
&\mathbb E\left [ \int_{\mathbb R^p}\left[ \hat g_N(\btheta_0) -g(\btheta_0|\bs_0)\right]^2\emph{d}\btheta_0\right]\\
& \quad =\left(\Phi_1(\bs_0)h_N^4+\Phi_2(\bs_0)\Delta(k_N)+\Phi_3(\bs_0)h_N^2D_4(k_N)+\frac{\int_{\mathbb R^p} K^2(\btheta)\emph{d}\btheta}{k_Nh_N^{p}}\right)\\
& \qquad \times\left (1+\emph{o}(1)\right).
\end{align*}
\item {\bf For $\boldsymbol{m\notin \{2,4\}}$},
\begin{align*}
&\mathbb E\left [ \int_{\mathbb R^p}\left[ \hat g_N(\btheta_0) -g(\btheta_0|\bs_0)\right]^2\emph{d}\btheta_0\right]\\
& \quad =\left(\Phi_1(\bs_0)h_N^4+\Phi_2(\bs_0)\Delta_m(k_N)+\Phi_3(\bs_0)h_N^2D_m(k_N)+\frac{\int_{\mathbb R^p} K^2(\btheta)\emph{d}\btheta}{k_Nh_N^{p}}\right)\\
& \qquad  \times\left (1+\emph{o}(1)\right).
\end{align*}
\end{enumerate}
 \end{theo}
By balancing the terms in Theorem \ref{theorem4}, we are led to the following useful corollary:
\begin{cor}[Rates of convergence]
\label{cor}
Under the conditions of Theorem \ref{theorem4}, one has:
\begin{enumerate}
\item {\bf For $\boldsymbol{m\in \{1,2,3\}}$}, there exists a sequence $\{k_N\}$ with $k_N\propto N^{\frac{p+4}{p+8}}$ and a sequence $\{h_N\}$ with $h_N\propto N^{-\frac{1}{p+8}}$ such that
\begin{align*}
&\mathbb E\left [ \int_{\mathbb R^p}\left[ \hat g_N(\btheta_0) -g(\btheta_0|\bs_0)\right]^2\emph{d}\btheta_0\right]\\
& \quad =\left(\frac{L^{4-m}\Phi_1(\bs_0)}{\xi_0(1-m/4)}+\Phi_2(\bs_0)+\int_{\mathbb R^p} K^2(\btheta)\emph{d}\btheta\right) N^{-\frac{4}{p+8}}+\emph{o}\left (N^{-\frac{4}{p+8}}\right).
\end{align*}
\item {\bf For $\boldsymbol{m=4}$}, there exists a sequence $\{k_N\}$ with $k_N\propto N^{\frac{p+4}{p+8}}$ and a sequence $\{h_N\}$ with $h_N\propto N^{-\frac{1}{p+8}}$ such that
\begin{align*}
&\mathbb E\left [ \int_{\mathbb R^p}\left[ \hat g_N(\btheta_0) -g(\btheta_0|\bs_0)\right]^2\emph{d}\btheta_0\right] \\
& \quad = \frac{4\Phi_1(\bs_0)}{\xi_0(p+8)} N^{-\frac{4}{p+8}}\log N+\emph{o}\left(N^{-\frac{4}{p+8}}\log N\right).
\end{align*}
\item {\bf For $\boldsymbol{m>4}$}, there exists a sequence $\{k_N\}$ with $k_N\propto N^{\frac{p+4}{m+p+4}}$ and a sequence $\{h_N\}$ with $h_N\propto N^{-\frac{1}{m+p+4}}$ such that
\begin{align*}
&\mathbb E\left [ \int_{\mathbb R^p}\left[ \hat g_N(\btheta_0) -g(\btheta_0|\bs_0)\right]^2\emph{d}\btheta_0\right]\\
& \quad =\left(\frac{m\Phi_1(\bs_0)}{\xi_0^{4/m}(m-4)}+\Phi_2(\bs_0)+\frac{m\Phi_3(\bs_0)}{\xi_0^{2/m}(m-2)}+\int_{\mathbb R^p} K^2(\btheta)\emph{d}\btheta\right)N^{-\frac{4}{m+p+4}}\\
& \qquad +\emph{o}\left(N^{-\frac{4}{m+p+4}}\right).
\end{align*}
\end{enumerate}
 \end{cor}
Several remarks are in order:
\begin{enumerate}
\item From a practical perspective, the fundamental problem is that of the joint choice of $k_N$ and $h_N$ in the absence of {\it a priori} information regarding the posterior $g(.|\bs_0)$. Various bandwidth selection
rules for conditional density estimates have been proposed in the literature \citep[see, e.g.,][]{Bashta,Hall,Fan}. However most if not all of these procedures pertain to kernel-type estimates and are difficult to adapt to our nearest-neighbor setting. Moreover, they are tailored to global statistical performance criteria, whereas the problem we are facing is local since $\bs_0$ is held fixed. Devising a good methodology to automatically select both parameters $k_N$ and $h_N$ in function of $\bs_0$ necessitates a specific analysis, which we believe is beyond the scope of the present paper.
\item Nevertheless, Corollary \ref{cor} provides a useful insight into the proportion of simulated values which should be accepted by the algorithm. For example, for $m>4$, a rough rule of thumb is obtained by taking $k_N \approx N^{(p+4)/(m+p+4)}$, so that a fraction of about $k_N/N \approx  N^{-m/(m+p+4)}$ {\sc abc}-simulations should not be rejected.
\end{enumerate}
 \section{Proofs}
\subsection{Proof of Proposition \ref{theorem1}}
Denote by $(\tilde\bTheta_1,\tilde \bS_1), \hdots,(\tilde\bTheta_k,\tilde
\bS_k)$ i.i.d. random couples with common probability density 
\begin{equation*}
\frac{1}{C_{d_{(k+1)}}} {\mathbf 1}_{[ \|
  \bs-\bs_0\|\leq d_{(k+1)}]} f(\btheta,\bs),
  \end{equation*}
 where the normalizing constant $C_{d_{(k+1)}}$ is defined by
 $$C_{d_{(k+1)}}= \int_{\mathbb R^p}\int_{\mathcal B_m(\bs_0,d_{(k+1)})}f(\btheta,\bs)\mbox{d}{\btheta}\mbox{d}{\bs}.$$
 Note, since $\bs_0$ belongs by assumption to the support of $\bar f$, that the constant $C_{d_{(k+1)}}$ is positive. To prove the first statement of the theorem, it is enough to establish that, for any test functions $\Phi$ and $\varphi$, with $\Phi$ symmetric in its arguments, one has
\begin{align*}
&\mathbb E \left [\Phi\left ((\bTheta_{(1)},\bS_{(1)}), \hdots,(\bTheta_{(k)},\bS_{(k)})\right)\varphi(d_{(k+1)})\right] \\
&\quad =  \mathbb E\left [\Phi\left ((\tilde\bTheta_1,\tilde \bS_1), \hdots,(\tilde\bTheta_k,\tilde
\bS_k)\right)\varphi(d_{(k+1)})\right].
\end{align*}
This can be achieved by adapting the proof of Lemma A.1 in \citet{cg1} to this context. Details are omitted. 
 \medskip

To prove the second statement, it suffices to show that, for any test functions $\Phi$ and $\varphi$ (with $\Phi$ not necessarily symmetric), one has
\begin{align*}
&\mathbb E [\Phi\left ((\bTheta_{(1)},\bS_{(1)}), \hdots,(\bTheta_{(k)},\bS_{(k)})\right)\varphi(d_{(k+1)})] \\
& \quad =\mathbb E[\Phi\left ((\tilde\bTheta_{(1)},\tilde \bS_{(1)}), \hdots,(\tilde\bTheta_{(k)},\tilde
\bS_{(k)})\right)\varphi(d_{(k+1)})].
\end{align*}
The arguments of \citet{cg1} may be repeated {\it mutatis mutandis} by replacing the $k$-com\-bi\-na\-tions of $\{1,\hdots,N\}$ by the $k$-permutations.
\subsection{Proof of Theorem \ref{theorem2}}
The proof strongly relies on Proposition \ref{theorem1}. It is assumed throughout that $\bs_0$ is a Lebesgue point of $\bar f$ ($\lambda_m$-almost all points satisfy this requirement) such that $\bar f(\bs_0)>0$. We note that this forces $\bs_0$ to belong to the support of $\bar f$, so that the assumption of Proposition \ref{theorem1} is satisfied. The collection of valid $\bs_0$ will vary during the proof, but only on subsets of Lebesgue measure 0. Similarly, we fix $\btheta_0 \in \mathbb R^p$, up to subsets of Lebesgue measure 0 which will appear in the proof.
\medskip

First observe that, according to Proposition \ref{theorem1}, 
$$\mathbb E[\hat g_N(\btheta_0)\, \big |\,d_{(k_N+1)}]=\frac{1}{C_{d_{(k_N+1)}}} \int_{\mathbb R^p} K_{h_N} (\btheta_0-\btheta)\left(\int_{\mathcal B_m(\bs_0,d_{(k_N+1)})}
f(\btheta,\bs)\mbox{d}\bs \right)\mbox{d}\btheta,$$
where, for any $\delta>0$, $C_{\delta}=\int_{\mathcal B_m(\bs_0,{\delta})}\bar f(\bs)\mbox{d}{\bs}$. Put differently, by Fubini's theorem,
\begin{equation}
\label{EC}
\mathbb E[\hat g_N(\btheta_0)\, \big |\,d_{(k_N+1)}] =\frac{1}{C_{d_{(k_N+1)}}} \int_{\mathbb R^p}\int_{\mathcal B_m(\bs_0,d_{(k_N+1)})} K_{h_N} (\btheta_0-\btheta)f(\btheta,\bs)\mbox{d}\btheta \mbox{d}\bs.
\end{equation}

The proof starts with the variance-bias decomposition 
\begin{align}
\mathbb E \left [\hat g_N(\btheta_0)-g(\btheta_0)\right]^2 &=  \mathbb E \left[ \mathbb E\left [\left (\hat g_N(\btheta_0) - \mathbb E[\hat
g_N(\btheta_0)\, \big | \,d_{(k_N+1)}]\right)^2\,\big | \, d_{(k_N+1)}\right]\right] \nonumber\\
& \quad +  \mathbb E\left [\mathbb E[\hat
g_N(\btheta_0)\,\big | \,d_{(k_N+1)}]-g(\btheta_0)\right]^2.\label{BV}
\end{align}
Our goal is to show that, under our assumptions, both terms on the right-hand side of (\ref{BV}) tend to 0 as $N \to \infty$. We start with the analysis of the second one, by noting that
\begin{align*}
&  \left |\mathbb E[\hat g_N(\btheta_0)\, \big |\,d_{(k_N+1)}]-g(\btheta_0)\right|  \\  
&\quad = \left |\frac{1}{C_{d_{(k_N+1)}}} \int_{\mathbb R^p}\int_{\mathcal B_m(\bs_0,d_{(k_N+1)})} K_{h_N} (\btheta_0-\btheta)f(\btheta,\bs)\mbox{d}\btheta\mbox{d}\bs -\frac{f(\btheta_0,\bs_0)}{\bar f(\bs_0)}\right|,
\end{align*}
where we used (\ref{EC}) and the definition of $g(\btheta_0)$. Equivalently,
\begin{align*}
&  \left |\mathbb E[\hat g_N(\btheta_0)\, \big |\,d_{(k_N+1)}]-g(\btheta_0)\right|  \\  
&\quad = \bigg |\frac{V_m\,d_{(k_N+1)}^m}{C_{d_{(k_N+1)}}} \frac{1}{V_m\,d_{(k_N+1)}^m}\int_{\mathbb R^p}\int_{\mathcal B_m(\bs_0,d_{(k_N+1)})} K_{h_N} (\btheta_0-\btheta)f(\btheta,\bs)\mbox{d}\btheta\mbox{d}\bs \\
& \quad \qquad-\frac{f(\btheta_0,\bs_0)}{\bar f(\bs_0)}\bigg|.
\end{align*}
For a fixed pair $(\btheta_0,\bs_0)$ and all $h,\delta>0$, set 
$$\psi_{\btheta_0,\bs_0}(h,\delta)=\left |\frac{V_m\delta^m}{C_{\delta}}\frac{1}{V_m\delta^m} \int_{\mathbb R^p}\int_{\mathcal B_m(\bs_0,\delta)} K_{h} (\btheta_0-\btheta)f(\btheta,\bs)\mbox{d}\btheta\mbox{d}\bs-\frac{f(\btheta_0,\bs_0)}{\bar f(\bs_0)}\right|.
$$
According to technical Lemma \ref{lemmetechnique1} $(i)$, the quantity ${V_m\delta^m}/{C_{\delta}}$ tends to $1/\bar f(\bs_0)$ as $\delta \to 0$. Therefore, by the first statement of Theorem \ref{WZ2}, we deduce that for $\lambda_p\otimes \lambda_m$-almost all pairs $(\btheta_0,\bs_0)\in \mathbb R^p\times \mathbb R^m$, $\lim_{(h,\delta)\to (0,0)}\psi^2_{\btheta_0,\bs_0}(h,\delta)= 0$.
\medskip

Next, introduce $\pi^{\star}$ (respectively, $f^{\star}$), the maximal function defined in Theorem \ref{WZ1} (respectively, Theorem \ref{WZ2}). Take any $\delta_0>0$. On the one hand, by the very definition of $f^{\star}$, 
$$\sup_{h>0,\delta_0\geq \delta>0}\left [\psi_{\btheta_0,\bs_0}(h,\delta)\right]\leq \sup_{0<\delta\leq \delta_0}\left [\frac{V_m\delta^m}{C_{\delta}}\right]f^{\star}(\btheta_0,\bs_0)+\frac{f(\btheta_0,\bs_0)}{\bar f(\bs_0)}.$$
On the other hand, for $\delta>\delta_0$,
$$\psi_{\btheta_0,\bs_0}(h,\delta) \leq \frac{1}{C_{\delta_0}}\int_{\mathbb R^p}K_{h} (\btheta_0-\btheta)\pi(\btheta)\mbox{d}\btheta+\frac{f(\btheta_0,\bs_0)}{\bar f(\bs_0)},$$
so that
$$\sup_{h>0,\delta>\delta_0}\left [\psi_{\btheta_0,\bs_0}(h,\delta)\right] \leq \frac{\pi^{\star}(\btheta_0)}{C_{\delta_0}}+\frac{f(\btheta_0,\bs_0)}{\bar f(\bs_0)}.$$
\medskip

Thus, putting all the pieces together, we infer that for $\lambda_p\otimes \lambda_m$-almost all pairs $(\btheta_0,\bs_0)\in \mathbb R^p\times \mathbb R^m$,
\begin{equation}
\label{SUP1}
\sup_{h>0,\delta>0}\left [\psi_{\btheta_0,\bs_0}(h,\delta)\right]\leq \sup_{0<\delta\leq \delta_0}\left [\frac{V_m\delta^m}{C_{\delta}}\right]f^{\star}(\btheta_0,\bs_0)+\frac{\pi^{\star}(\btheta_0)}{C_{\delta_0}}+\frac{2f(\btheta_0,\bs_0)}{\bar f(\bs_0)}.
\end{equation}
In consequence, by Lemma \ref{lemmetechnique1} $(ii)$, Theorem \ref{WZ1} $(ii)$ and Theorem \ref{WZ2} $(ii)$, for such pairs $(\btheta_0,\bs_0)$,
\begin{equation}
\label{SUP2}
\sup_{h>0,\delta>0}\left [\psi^2_{\btheta_0,\bs_0}(h,\delta)\right]<\infty.
\end{equation}

Now, since $d_{(k_N+1)}\to 0$ with probability 1 whenever $k_N/N\to 0$ \citep[see, e.g., Lemma 5.1 in][]{DGL}, 
we conclude by Le\-bes\-gue's dominated convergence theorem that the bias term in (\ref{BV}) tends to 0 as $N\to \infty$.
\medskip

To finish the proof, it remains to show that the first term of
(\ref{BV}) vanishes as $N\to \infty$. This is easier. Just note that, using again Proposition \ref{theorem1},
\begin{align}\label{variance}
&\mathbb E\left [\left (\hat g_N(\btheta_0) - \mathbb E[\hat
g_N(\btheta_0)\, \big | \,d_{(k_N+1)}]\right)^2\,\big | \, d_{(k_N+1)}\right]\nonumber\\
& \quad = \frac{1}{k_Nh_N^{2p}}\frac{1}{C_{d_{(k_N+1)}}} \int_{\mathbb R^p}K^2\left(\frac{\btheta_0-\btheta}{h_N}\right) \left(\int_{\mathcal B_m(\bs_0,d_{(k_N+1)})}f(\btheta,\bs)\mbox{d}\bs \right)\mbox{d}\btheta\nonumber\\
& \quad \quad -\frac{1}{k_N}\left(\mathbb E \left [\hat g_N(\btheta_0)\,\big | \, d_{(k_N+1)}\right]\right)^2.
\end{align}
Hence, if $K$ is bounded by, say, $\|K\|_{\infty}$,
\begin{align*}
&\mathbb E\left [\left (\hat g_N(\btheta_0) - \mathbb E[\hat
g_N(\btheta_0)\, \big | \,d_{(k_N+1)}]\right)^2\,\big | \, d_{(k_N+1)}\right]\\
& \quad \leq \frac{1}{k_Nh_N^{2p}}\frac{1}{C_{d_{(k_N+1)}}} \int_{\mathbb R^p}K^2\left(\frac{\btheta_0-\btheta}{h_N}\right) \left(\int_{\mathcal B_m(\bs_0,d_{(k_N+1)})}f(\btheta,\bs)\mbox{d}\bs \right)\mbox{d}\btheta\\
& \quad \leq\frac{1}{k_Nh_N^{p}}\frac{\|K\|_\infty}{C_{d_{(k_N+1)}}} \int_{\mathbb
  R^p}\int_{\mathcal B_m(\bs_0,d_{(k_N+1)})}K_{h_N}(\btheta_0-\btheta)
f(\btheta,\bs)\mbox{d}\btheta \mbox{d}\bs.
\end{align*}
Thus, using (\ref{SUP2}), we obtain
$$\mathbb E\left [\left (\hat g_N(\btheta_0) - \mathbb E[\hat
g_N(\btheta_0)\, \big | \,d_{(k_N+1)}]\right)^2\,\big | \, d_{(k_N+1)}\right]\leq \frac{ C}{k_Nh_N^p}$$
for some positive constant $C$ depending on $\btheta_0$, $\bs_0$ and
$K$, but independent of $h_N$ and $k_N$. This shows that the variance term goes to 0 as $k_Nh_N^p \to \infty$ and concludes the proof of the theorem.
\subsection{Proof of Theorem \ref{theorem3}}
We start as in the proof of Theorem \ref{theorem2} and write, using Fubini's theorem,
\begin{align}
&\mathbb E\left[\int_{\mathbb R^p}\left[\hat g_N(\btheta_0)-g(\btheta_0)\right]^2\mbox{d}\btheta_0\right]\nonumber\\
& \quad =\mathbb E\left[\int_{\mathbb R^p}\mathbb E\left [\left (\hat g_N(\btheta_0) - \mathbb E[\hat
g_N(\btheta_0)\, \big | \,d_{(k_N+1)}]\right)^2\,\big | \, d_{(k_N+1)}\right]\mbox{d}\btheta_0\right]\nonumber\\
& \qquad + \mathbb E\left[\int_{\mathbb R^p} \left[\mathbb E[\hat
g_N(\btheta_0)\,\big | \,d_{(k_N+1)}]-g(\btheta_0)\right]^2\mbox{d}\btheta_0\right].\label{BVI}
\end{align}
It has already been seen that
\begin{align*}
&\mathbb E\left [\left (\hat g_N(\btheta_0) - \mathbb E[\hat
g_N(\btheta_0)\, \big | \,d_{(k_N+1)}]\right)^2\,\big | \, d_{(k_N+1)}\right]\\
& \quad \leq \frac{1}{k_Nh_N^{2p}}\frac{1}{C_{d_{(k_N+1)}}} \int_{\mathbb R^p}\int_{\mathcal B_m(\bs_0,d_{(k_N+1)})}K^2\left(\frac{\btheta_0-\btheta}{h_N}\right) f(\btheta,\bs)\mbox{d}\btheta \mbox{d}\bs.
\end{align*}
Consequently, by definition of $C_{d_{(k_N+1)}}$, we are led to
\begin{equation*}
\int_{\mathbb R^p}\mathbb E\left [\left (\hat g_N(\btheta_0) - \mathbb E[\hat
g_N(\btheta_0)\, \big | \,d_{(k_N+1)}]\right)^2\,\big | \, d_{(k_N+1)}\right]\mbox{d}\btheta_0 \leq \frac{\int_{\mathbb R^p} K^2(\btheta)\mbox{d}\btheta}{k_Nh_N^{p}}.
\end{equation*}
This shows that the first term in (\ref{BVI}) tends to 0 as $k_Nh_N^p\to \infty$. 
\medskip

Let us now turn to the analysis of the bias term. With the notation of the proof of Theorem \ref{theorem2}, we may write
 $$\mathbb E \left [\int_{\mathbb R^p}\left[\mathbb E[\hat g_N(\btheta_0)\, \big |\,d_{(k_N+1)}]-g(\btheta_0)\right]^2\mbox{d}\btheta_0 \right] = \mathbb E \left [ \int_{\mathbb R^p} \psi^2_{\btheta_0,\bs_0}(h_N,d_{(k_N+1)})\mbox{d}\btheta_0\right].$$
It is known from the proof of Theorem \ref{theorem2} that the limit of $\psi^2_{\btheta_0,\bs_0}(h,\delta)$ is 0 for $\lambda_p\otimes \lambda_m$-almost all $(\btheta_0,\bs_0)\in \mathbb R^p\times \mathbb R^m$, whenever $(h,\delta)\to (0,0)$.
Take any $\delta_0>0$. Denoting by $f^{\star}$ (respectively, $\pi^{\star}$) the maximal function defined in Theorem \ref{WZ2} (respectively, Theorem \ref{WZ1}), we also know (inequality (\ref{SUP1})) that 
$$
\sup_{h>0,\delta>0}\left [\psi_{\btheta_0,\bs_0}(h,\delta)\right]\leq \sup_{0<\delta\leq \delta_0}\left [\frac{V_m\delta^m}{C_{\delta}}\right]f^{\star}(\btheta_0,\bs_0)+\frac{\pi^{\star}(\btheta_0)}{C_{\delta_0}}+\frac{2f(\btheta_0,\bs_0)}{\bar f(\bs_0)}.
$$
Thus, because $(a+b+c)^2\leq 3 a^2 +3b^2+3c^2$,
\begin{align*}
&\sup_{h>0,\delta>0}\left [\psi^2_{\btheta_0,\bs_0}(h,\delta)\right]\\
& \quad \leq 3\left(\sup_{0<\delta \leq \delta_0}\left [\frac{V_m\delta^m}{C_{\delta}}\right]f^{\star}(\btheta_0,\bs_0)\right)^2+3\left(\frac{\pi^{\star}(\btheta_0)}{C_{\delta_0}}\right)^2+12\left(\frac{f(\btheta_0,\bs_0)}{\bar f(\bs_0)}\right)^2.
\end{align*}
By Lemma \ref{lemmetechnique1} $(ii)$, the supremum on the right-hand side is bounded. Moreover, by assumption, $f$ is in $L^2(\mathbb R^p \times \mathbb R^m)$. Therefore the function $\btheta_0 \mapsto f(\btheta_0,\bs_0)$ is in $L^2(\mathbb R^p)$ as well for $\lambda_m$-almost all $\bs_0 \in \mathbb R^m$. Similarly, for $\lambda_m$-almost all $\bs_0$, by Theorem \ref{WZ2} $(iii)$, the function $\btheta_0 \mapsto f^{\star}(\btheta_0,\bs_0)$ is in $L^2(\mathbb R^p)$. Finally, $\pi^{\star}$ belongs to $L^2(\mathbb R^p)$ by Theorem \ref{WZ1} $(iii)$.
Since $d_{(k_N+1)}\to 0$ with probability 1 whenever $k_N/N\to 0$, the conclusion follows from Lebesgue's dominated convergence theorem.
 \subsection{Proof of Theorem \ref{theorem4}}
Throughout the proof, it is assumed that the Lebesgue point $\bs_0$ is fixed and such that $\bar f(\bs_0)>0$. This forces $\bs_0$ to belong to the support of $\bar f$.
\medskip

As in the proofs of Theorem \ref{theorem2} and Theorem \ref{theorem3}, we set, for any $\btheta_0 \in \mathbb R^p$ and all $h,\delta>0$,
$$\psi_{\btheta_0,\bs_0}(h,\delta)=\left |\frac{V_m\delta^m}{C_{\delta}}\frac{1}{V_m\delta^m} \int_{\mathbb R^p}\int_{\mathcal B_m(\bs_0,\delta)} K_{h} (\btheta_0-\btheta)f(\btheta,\bs)\mbox{d}\btheta\mbox{d}\bs-\frac{f(\btheta_0,\bs_0)}{\bar f(\bs_0)}\right|,
$$
where $C_{\delta}=\int_{\mathcal B_m(\bs_0,{\delta})}\bar f(\bs)\mbox{d}{\bs}$. With this notation, it is readily seen from identities (\ref{variance}) and (\ref{BVI}) that
\begin{align*}
&\mathbb E\left [ \int_{\mathbb R^p}\left[ \hat g_N(\btheta_0) -g(\btheta_0|\bs_0)\right]^2\mbox{d}\btheta_0\right] \\
&\quad =\mathbb E \left [ \int_{\mathbb R^p} \psi^2_{\btheta_0,\bs_0}(h_N,d_{(k_N+1)})\mbox{d}\btheta_0\right]+\frac{\int_{\mathbb R^p} K^2(\btheta)\mbox{d}\btheta}{k_Nh_N^{p}}\\
& \quad\quad  -\frac{1}{k_N}\mathbb E \left [ \int_{\mathbb R^p} \left(\mathbb E \left [\hat g_N(\btheta_0)\,\big | \, d_{(k_N+1)}\right]\right)^2 \mbox{d}\btheta_0\right].
\end{align*}
Recall that
$$\mathbb E[\hat g_N(\btheta_0)\, \big |\,d_{(k_N+1)}]=\frac{1}{C_{d_{(k_N+1)}}} \int_{\mathbb R^p} K_{h_N} (\btheta_0-\btheta)\left(\int_{\mathcal B_m(\bs_0,d_{(k_N+1)})}
f(\btheta,\bs)\mbox{d}\bs \right)\mbox{d}\btheta,$$
and the same arguments as in the proof of Theorem \ref{theorem3} reveal that
$$\sup_{h_N>0,L\geq d_{(k_N+1)}>0}\left(\mathbb E \left [\hat g_N(\btheta_0)\,\big |\,d_{(k_N+1)}\right]\right)^2\leq \left(\sup_{0<\delta \leq L}\left [\frac{V_m\delta^m}{C_{\delta}}\right]f^{\star}(\btheta_0,\bs_0)\right)^2.$$
Since $f$ is in $L^2(\mathbb R^p\times \mathbb R^m)$ by Assumption $[\mathbf{A2}]$, this ensures that for $\lambda_m$-almost all $\bs_0\in \mathbb R^m$,
$$\mathbb E \left [ \int_{\mathbb R^p} \left(\mathbb E \left [\hat g_N(\btheta_0)\,\big | \, d_{(k_N+1)}\right]\right)^2 \mbox{d}\btheta_0\right]<\infty$$
and 
$$\frac{1}{k_N}\mathbb E \left [ \int_{\mathbb R^p} \left(\mathbb E \left [\hat g_N(\btheta_0)\,\big | \, d_{(k_N+1)}\right]\right)^2 \mbox{d}\btheta_0\right]=\mbox{O}\left(\frac{1}{k_N}\right).$$
In particular, 
$$\frac{1}{k_N}\mathbb E \left [ \int_{\mathbb R^p} \left(\mathbb E \left [\hat g_N(\btheta_0)\,\big | \, d_{(k_N+1)}\right]\right)^2 \mbox{d}\btheta_0\right]=\mbox{o}\left(\frac{1}{k_Nh_N^{p}}\right).$$
The rest of the proof is devoted to the study of the rate of convergence to 0 of the quantity
$$\mathbb E \left [ \int_{\mathbb R^p} \psi^2_{\btheta_0,\bs_0}(h_N,d_{(k_N+1)})\mbox{d}\btheta_0\right].$$
By an elementary change of variables, using the symmetry of $K$,
\begin{align*}
&\frac{1}{V_m\delta^m}\int_{\mathbb R^p}\int_{\mathcal B_m(\bs_0,\delta)} K_{h} (\btheta_0-\btheta)f(\btheta,\bs)\mbox{d}\btheta\mbox{d}\bs\\
&\quad = \frac{1}{V_m}\int_{\mathbb R^p}\int_{\mathcal B_m(\mathbf 0,1)} K (\btheta)f(\btheta_0+h\btheta,\bs_0+\delta\bs)\mbox{d}\btheta\mbox{d}\bs.
\end{align*}
Next, by the multivariate Taylor's theorem applied to $f$ around $(\btheta_0,\bs_0)$ (which is valid here by Assumption $[\mathbf{A3}]$),
\begin{align*}
f(\btheta_0+h\btheta,\bs_0+\delta\bs)=&f(\btheta_0,\bs_0)+\sum_{|\beta|=1}D^{\beta}f(\btheta_0,\bs_0)(h\btheta,\delta\bs)^\beta\\
& +\sum_{|\beta|=2}\frac{D^{\beta}f(\btheta_0,\bs_0)}{{\beta !}}(h\btheta,\delta\bs)^\beta\\
&+\sum_{|\beta|=3}R_{\beta}(\btheta_0+h\btheta,\bs_0+\delta\bs)(h\btheta,\delta\bs)^\beta,
\end{align*}
where each component of the remainder term takes the form
$$R_{\beta}(\btheta_0+h\btheta,\bs_0+\delta\bs)= \frac{3}{\beta!}\int_0^1 (1-t)^2D^\beta f(\btheta_0+th\btheta,\bs_0+t\delta\bs)\mbox{d}t.$$
In view of the symmetry of $K$ and the ball $\mathcal B_m(\mathbf 0,1)$, it is clear that
$$\int_{\mathbb R^p}\int_{\mathcal B_m(\mathbf 0,1)} K (\btheta) \sum_{|\beta|=1}D^{\beta}f(\btheta_0,\bs_0)(h\btheta,\delta\bs)^{\beta} \mbox{d}\btheta\mbox{d}\bs=0.$$
Similarly, elementary calculations reveal that
\begin{align*}
&\frac{1}{V_m}\int_{\mathbb R^p}\int_{\mathcal B_m(\mathbf 0,1)} K (\btheta) \sum_{|\beta|=2}\frac{D^{\beta}f(\btheta_0,\bs_0)}{\beta!}(h\btheta,\delta\bs)^\beta\mbox{d}\btheta\mbox{d}\bs\\
&\quad = \phi_1(\btheta_0,\bs_0)h^2+\phi_2(\btheta_0,\bs_0)\delta^2
\end{align*}
(where $\phi_1$ is defined in the statement of Theorem \ref{theorem4}), and 
\begin{align*}
\phi_2(\btheta_0,\bs_0) &=\frac{1}{2V_m}\sum_{j=1}^m\frac{\partial^2 f(\btheta_0,\bs_0)}{\partial s_j^2}\int_{\mathcal B_m(\mathbf 0,1)} s_j^2\mbox{d}\bs.
\end{align*}
Using expression (\ref{volume}) of $V_m$, an elementary verification shows that
$$\frac{1}{V_m}\int_{\mathcal B_m(\mathbf 0,1)} s_j^2\mbox{d}\bs=\frac{1}{m+2}\quad \mbox{and}\quad \phi_2(\btheta_0,\bs_0)=\frac{1}{2m+4}\sum_{j=1}^m\frac{\partial^2 f(\btheta_0,\bs_0)}{\partial s_j^2}.$$
Let us now define $(\mathbf h,\bd{\delta})=(h,\dots,h,\delta,\dots,\delta)$ (where $h$ is replicated $p$ times and $\delta$ is replicated $m$ times) and care about the remainder term $R_{\beta}(\btheta_0+h\btheta,\bs_0+\delta\bs)$. For any multi-index $\beta$ with $|\beta|=3$,  it holds
$$\int_{\mathbb R^p}\int_{\mathcal B_m(\mathbf 0,1)} K (\btheta)R_{\beta}(\btheta_0+h\btheta,\bs_0+\delta\bs)(h\btheta,\delta\bs)^\beta\mbox{d}\btheta\mbox{d}\bs=(\mathbf h,\bd{\delta})^\beta A_\beta(\btheta_0,h,\delta),$$
where, by definition,
$$ A_\beta(\btheta_0,h,\delta)=\int_{\mathbb R^p}\int_{\mathcal B_m(\mathbf 0,1)} K (\btheta)R_{\beta}(\btheta_0+h\btheta,\bs_0+\delta\bs)(\btheta,\bs)^\beta\mbox{d}\btheta\mbox{d}\bs.$$
[Note that $A_\beta(\btheta_0,h,\delta)$ depends in fact upon $\bs_0$ as well, but since this dependency is not crucial, we leave it out in the notation.]
Finally, 
\begin{align*}
&\frac{1}{V_m\delta^m}\int_{\mathbb R^p}\int_{\mathcal B_m(\bs_0,\delta)} K_{h} (\btheta_0-\btheta)f(\btheta,\bs)\mbox{d}\btheta\mbox{d}\bs\\
& \quad = f(\btheta_0,\bs_0)+ \phi_1(\btheta_0,\bs_0)h^2+ \phi_2(\btheta_0,\bs_0)\delta^2+ \sum_{|\beta|=3}(\mathbf h,\bd{\delta})^\beta A_\beta(\btheta_0,h,\delta).
\end{align*}
Considering now the function
$$\tau_{\bs_0}(\delta)=\frac{C_{\delta}}{V_m\delta^m}=\frac{1}{V_m\delta^m}\int_{\mathcal B_m(\bs_0,{\delta})}\bar f(\bs)\mbox{d}{\bs}=\frac{1}{V_m}\int_{\mathcal B_m(\mathbf 0,1)}\bar f(\bs_0+\delta\bs)\mbox{d}{\bs},$$
and the asymptotic expansion of $1/\tau_{\bs_0}$ around 0, a similar analysis shows that
$$\frac{V_m\delta^m}{C_{\delta}}=\frac{1}{\bar f(\bs_0)}-\frac{\phi_3(\bs_0)}{{\bar f}^2(\bs_0)}\delta^2+ \delta^3 \zeta_1(\delta)$$
(where $\phi_3$ is defined in the statement of Theorem \ref{theorem4}), and, with a slight abuse of notation, there exists $t\in(0,1)$ such that $\zeta_1(\delta)=H(t\delta)/\tau_{\bs_0}^4(t\delta)$. In this last expression, the function $H$ depends only on the successive derivatives $D^\beta\bar f(\bs_0+t\delta\bs)$ for $0\leq|\beta|\leq 3$ and is therefore bounded thanks to Assumption $[{\bf A1}]$. Besides, by the very definition of $\xi_0$ and technical Lemma \ref{min}, 
$$\tau_{\bs_0}(t\delta)=\frac{1}{V_m(t\delta)^m}\int_{\mathcal B_m(\bs_0,{t\delta})}\bar f(\bs)\mbox{d}{\bs}\geq \frac{\xi_0}{V_m}>0.$$
Thus, the function $\zeta_1(\delta)$ is such that $\sup_{0<\delta\leq L}\zeta_1(\delta)<\infty$. Putting all the pieces together, we conclude that
$$\psi_{\btheta_0,\bs_0}(h,\delta)=\left|\phi_4(\btheta_0,\bs_0) h^2+\phi_5(\btheta_0,\bs_0)\delta^2+ h^2 \zeta_2(\btheta_0,h,\delta)+\delta^2\zeta_3(\btheta_0,h,\delta)\right|,$$
where 
$$\phi_4(\btheta_0,\bs_0)=\frac{\phi_1(\btheta_0,\bs_0)}{\bar f(\bs_0)} \quad \mbox{and} \quad \phi_5(\btheta_0,\bs_0)=\frac{\phi_2(\btheta_0,\bs_0)\bar f(\bs_0)-\phi_3(\bs_0)f(\btheta_0,\bs_0)}{{\bar f}^2(\bs_0)}.$$
Moreover, one can check, using Assumption $[\mathbf{A2}]$ and the second statement of Assumption $[{\bf A3}]$ together with technical Lemma \ref{chiant}, that for $i=2,3$, $\zeta_i(\btheta_0,h,\delta)\to 0$ as $(h,\delta)\to (0,0)$, and 
$$\sup_{0<h<M, 0<\delta\leq L}\int_{\mathbb R^p}\zeta_i^2(\btheta_0,h,\delta)\mbox{d}\btheta_0<\infty$$
for all positive $M$. As a consequence, 
$$\int_{\mathbb R^p} \psi^2_{\btheta_0,\bs_0}(h,\delta)\mbox{d}\btheta_0= \Phi_1(\bs_0)h^4+\Phi_2(\bs_0)\delta^4+\Phi_3(\bs_0)h^2\delta^2+(h^2+\delta^2)^2\zeta_4(h,\delta)$$
($\Phi_1$, $\Phi_2$ and $\Phi_3$ are defined in the statement of Theorem \ref{theorem4}). Besides, for all positive $M$,
\begin{equation}
\label{RLB}
\sup_{0<h<M,0<\delta\leq L}\zeta_4(h,\delta)<\infty\quad\mbox{and } \lim_{(h,\delta)\to (0,0)}\zeta_4(h,\delta)= 0.
\end{equation}
Finally,
\begin{align*}
&\mathbb E\left [ \int_{\mathbb R^p}\left[ \hat g_N(\btheta_0) -g(\btheta_0|\bs_0)\right]^2\mbox{d}\btheta_0\right]\\
& \quad = \Phi_1(\bs_0)h_N^4+\Phi_2(\bs_0)\mathbb E[d_{(k_N+1)}^4]+\Phi_3(\bs_0)h_N^2\mathbb  E[d^2_{(k_N+1)}]+\frac{\int_{\mathbb R^p} K^2(\btheta)\mbox{d}\btheta}{k_Nh_N^{p}}\\
& \qquad +\mathbb E\left [(h_N^2+d_{(k_N+1)}^2)^2\zeta_4(h_N,d_{(k_N+1)})\right]+\mbox{o}\left(\frac{1}{k_Nh_N^{p}}\right).
\end{align*}
The conclusion is then an immediate consequence of (\ref{RLB}) and Assumption $[\mathbf{A1}$], together with Proposition \ref{distances} and Proposition \ref{distancesbis}, which respectively provide upper bounds on $\mathbb E[d_{(k_N+1)}^2]$ and $\mathbb E[d_{(k_N+1)}^4]$ depending on the dimension $m$.
\section{Some technical results}
\begin{lem}
\label{lemmetechnique1} Let $\bs_0 \in \mathbb R^m$ be a Lebesgue point of $\bar f$ such that $\bar f(\bs_0)>0$. For any $\delta>0$, let $C_{\delta}= \int_{\mathcal B_m(\bs_0,\delta)}\bar f(\bs)\emph{d}{\bs}$. One has
\begin{enumerate}
\item[$(i)$] $\lim_{\delta \to 0}V_m\delta^m/C_{\delta} = 1/\bar f(\bs_0)$.
\item[$(ii)$] For any $\delta_0>0$, $\sup_{0<\delta\leq \delta_0} V_m\delta^m/C_{\delta}<\infty$.
\end{enumerate}
\end{lem}
\noindent{\bf Proof of Lemma \ref{lemmetechnique1}}\quad The first statement is an immediate consequence of Lebesgue's differentiation theorem \citep[][Theorem 7.2]{Wheeden}. Take now $\delta_0>0$. Since $\bar f(\bs_0)>0$, it is routine to verify that the mapping $\delta \mapsto\frac{V_m\delta^m}{C_{\delta}}$ is positive and continuous on $(0,\delta_0]$. Thus, by $(i$), we deduce that $\sup_{0<\delta \leq \delta_0} V_m\delta^m/C_{\delta}<\infty$.
\hfill $\blacksquare$
\begin{lem}
\label{chiant}
Assume that the joint probability density $f$  is three times continuously differentiable on $\mathbb R^p \times \mathbb R^m$, and let $\beta$ be a multi-index satisfying $|\beta|=3$. Assume that $\sup_{\bs\in\mathbb R^m} \int_{\mathbb R^p}\left[D^\beta f(\btheta,\bs)\right]^2\emph{d}\btheta<\infty$, and, for $h,\delta>0$, consider the parameterized mapping $\btheta_0\mapsto A_\beta(\btheta_0,h,\delta)$, where
$$A_\beta(\btheta_0,h,\delta)=\int_{\mathbb R^p}\int_{\mathcal B_m(\mathbf 0,1)} K (\btheta)R_{\beta}(\btheta_0+h\btheta,\bs_0+\delta\bs)(\btheta,\bs)^\beta\emph{d}\btheta\emph{d}\bs,$$
with
$$R_{\beta}(\btheta_0+h\btheta,\bs_0+\delta\bs)= \int_0^1 (1-t)D^\beta f(\btheta_0+th\btheta,\bs_0+t\delta\bs)\emph{d}t.$$
Then 
$$\sup_{h,\delta>0}\int_{\mathbb R^p}A^2_\beta(\btheta_0,h,\delta)\emph{d}\btheta_0<\infty.$$
\end{lem}
\noindent{\bf Proof of Lemma \ref{chiant}}\quad The proof relies on an application of the generalized Minkowski's inequality \citep[see, e.g., Theorem 202 in][]{hlp}. Indeed,
\begin{align*}
&\left(\int_{\mathbb R^p}A^2_\beta(\btheta_0,h,\delta)\mbox{d}\btheta_0\right)^{\frac{1}{2}}\\
&\quad\leq \int_{\mathbb R^p}\int_{\mathcal B_m({\bf 0},1)}\int_0^1{\Sigma^{1/2}_\beta(\btheta,\bs,t)}(1-t)K(\btheta)\left|(\btheta,\bs)^\beta\right| \mbox{d}\btheta\mbox{d}\bs\mbox{d}t,
\end{align*}
where
$$\Sigma_\beta(\btheta,\bs,t)=\int_{\mathbb R^p} \left[ D^\beta f(\btheta_0+th\btheta,\bs_0+t\delta\bs)\right]^2 \mbox{d}\btheta_0.$$
Letting $C^2=\sup_{\bs\in\mathbb R^m} \int_{\mathbb R^p}[D^\beta f(\btheta,\bs)]^2\mbox{d}\btheta<\infty$, we obtain
$$\left(\int_{\mathbb R^p}A^2_\beta(\btheta_0,h,\delta)\mbox{d}\btheta_0\right)^{\frac{1}{2}}\leq C\int_{\mathbb R^p}\int_{\mathcal B_m({\bf 0},1)}\int_0^1(1-t)K(\btheta)\left|(\btheta,\bs)^\beta\right|\mbox{d}\btheta\mbox{d}\bs \mbox{d}t.$$
This upper bound is finite thanks to Assumption $[\mathbf{K2}$], and independent of $h$ and $\delta$. 
 \hfill $\blacksquare$
\begin{lem}
\label{min}
Let $\bs_0$ be a Lebesgue point of $\bar f$ such that $\bar f(\bs_0)>0$. Then, for all positive $L$,
$$0<\inf_{0<\delta\leq L}\frac{1}{\delta^m} \int_{\mathcal B_m(\bs_0,\delta)}\bar f(\bs)\emph{d}{\bs}<\infty.$$
\end{lem}
\noindent{\bf Proof of Lemma \ref{min}}\quad By exploiting the fact that $\bs_0$ is a Lebesgue point of $\bar f$ satisfying $\bar f(\bs_0)>0$, we deduce that for some positive $\delta_0<L$,
$$0<\inf_{0<\delta\leq \delta_0}\frac{1}{\delta^m} \int_{\mathcal B_m(\bs_0,\delta)}\bar f(\bs)\mbox{d}{\bs}<\infty.$$
Moreover, 
$$\frac{1}{L^m} \int_{\mathcal B_m(\bs_0,\delta_0)}\bar f(\bs)\mbox{d}{\bs} \leq \inf_{\delta_0<\delta \leq L}\frac{1}{\delta^m} \int_{\mathcal B_m(\bs_0,\delta)}\bar f(\bs)\mbox{d}{\bs}\leq \frac{1}{\delta_0^m}.$$
The quantity on the left-hand side is positive since $\bs_0$ belongs to the support of ${\bar f}$. This concludes the proof.
 \hfill $\blacksquare$
\begin{pro}
\label{distances}
Assume that the support of ${\bar f}$ is compact with diameter $L>0$.
Let $\bs_0$ be a Lebesgue point of $\bar f$ such that $\bar f(\bs_0)>0$. Set
$$\xi_0=\inf_{0<\delta\leq L}\frac{1}{\delta^m} \int_{\mathcal B_m(\bs_0,\delta)}\bar f(\bs)\emph{d}{\bs}.$$
Whenever $\frac{k_N+1}{N+1}\leq \xi_0 L^m$, one has:
\begin{enumerate}
\item For $m=2$,
$$\mathbb E \left [d_{(k_N+1)}^2\right]\leq\frac{1}{\xi_0}\left(1+\log \left (\xi_0 L^2 \frac{N+1}{k_N+1}\right)\right)\frac{k_N+1}{N+1}.$$
\item For $m\neq 2$,
$$\mathbb E \left [d_{(k_N+1)}^2\right]\leq\frac{m}{\xi_0^{2/m}(m-2)}\left (\frac{k_N+1}{N+1}\right)^{2/m} -\frac{L^{2-m}}{\xi_0(m/2-1)}\frac{k_N+1}{N+1}.$$
\end{enumerate}
\end{pro}
\noindent{\bf Proof of Proposition \ref{distances}}\quad First note, according to Lemma \ref{min}, that $0<\xi_0<\infty$. Next, observe that
$$\mathbb E \left [d_{(k_N+1)}^2\right]=\int_0^{L^2} \mathbb P \left \{ d_{(k_N+1)}> \sqrt \delta \right\}\mbox{d}\delta.$$
For some fixed $a \in (0,L^2)$, we use the decomposition
\begin{align*}
&\int_0^{L^2} \mathbb P \left \{ d_{(k_N+1)}> \sqrt \delta \right\}\mbox{d}\delta\\
&\quad =\int_0^{a} \mathbb P \left \{ d_{(k_N+1)}> \sqrt \delta \right\}\mbox{d}\delta+ \int_a^{L^2} \mathbb P \left \{ d_{(k_N+1)}> \sqrt \delta \right\}\mbox{d}\delta\\
& \quad \leq a + \int_a^{L^2} \mathbb P \left \{ d_{(k_N+1)}> \sqrt \delta \right\}\mbox{d}\delta.
\end{align*}
Introduce $p_0(\sqrt \delta)=\int_{\mathcal B_m(\bs_0,\sqrt \delta)}\bar f(\bs)\mbox{d}{\bs}$, which is positive since $\bs_0$ is in the support of $\bar f$. Using a binomial argument, we see that
\begin{align*}
\mathbb P \left \{ d_{(k_N+1)}> \sqrt \delta \right\}&=\sum_{j=0}^{k_N} {{N} \choose {j}} \left[p_0(\sqrt \delta)\right]^j \left [1-p_0(\sqrt \delta)\right]^{N-j}\\
& = \frac{1}{p_0(\sqrt \delta)} \sum_{j=0}^{k_N} {{N} \choose {j}} \left[p_0(\sqrt \delta)\right]^{j+1} \left [1-p_0(\sqrt \delta)\right]^{N-j}.
\end{align*}
By applying Lemma 3.1 in \citet{bcg2}, we obtain 
$$\mathbb P \left \{ d_{(k_N+1)}> \sqrt \delta \right\} \leq \frac{k_N+1}{N+1} \times \frac{1}{p_0(\sqrt \delta)}.$$
Consequently,
$$
\mathbb E \left [d_{(k_N+1)}^2\right] \leq a + \frac{1}{\xi_0} \frac{k_N+1}{N+1}\int_a^{L^2} \delta^{-m/2}\mbox{d}{\delta}.
$$
The conclusion is easily obtained by optimizing the right-hand side with respect to the parameter $a$.
 \hfill $\blacksquare$
 \begin{pro}
\label{distancesbis}
Assume that the support of ${\bar f}$ is compact with diameter $L>0$.
Let $\bs_0$ be a Lebesgue point of $\bar f$ such that $\bar f(\bs_0)>0$. Set
$$\xi_0=\inf_{0<\delta\leq L}\frac{1}{\delta^m} \int_{\mathcal B_m(\bs_0,\delta)}\bar f(\bs)\emph{d}{\bs}.$$
Whenever $\frac{k_N+1}{N+1}\leq \xi_0 L^m$, one has:
\begin{enumerate}
\item For $m=4$,
$$\mathbb E \left [d_{(k_N+1)}^4\right]\leq\frac{1}{\xi_0}\left(1+\log \left (\xi_0 L^4 \frac{N+1}{k_N+1}\right)\right)\frac{k_N+1}{N+1}.$$
\item For $m\neq 4$,
$$\mathbb E \left [d_{(k_N+1)}^4\right]\leq\frac{m}{\xi_0^{4/m}(m-4)}\left (\frac{k_N+1}{N+1}\right)^{4/m} -\frac{L^{4-m}}{\xi_0(m/4-1)}\frac{k_N+1}{N+1}.$$
\end{enumerate}
\end{pro}
\noindent{\bf Proof of Proposition \ref{distancesbis}}\quad Proof is similar to the one of Proposition \ref{distances}, and is therefore omitted.
 \hfill $\blacksquare$
\appendix
\section{Complements on singular integrals}
Recall that the convolution \citep[][Chapter 6]{Wheeden} of two measurable functions $f$ and $g$ in $\mathbb R^n$ is defined by
$$(f \star g)(\bx)=\int_{\mathbb  R^n}f(\mathbf y)g(\bx-\mathbf y)\mbox{d}\mathbf y,\quad \bx \in \mathbb R^n,
$$
provided the integral exists. This appendix is devoted to the study of some properties of convolution when $\mathbb R^n=\mathbb R^{n_1}\times \mathbb R^{n_2}$ and $g$ is of the form 
$$\varphi_{\varepsilon_1,\varepsilon_2}(\bx)=\frac{1}{\varepsilon_1^{n_1}\varepsilon_2^{n_2}}
\varphi_1\left (\frac{\bx_1}{\varepsilon_1}\right)\varphi_2\left (\frac{\bx_2}{\varepsilon_2}\right), \quad \bx=(\bx_1, \bx_2) \in \mathbb R^{n_1}\times \mathbb R^{n_2}.$$
More precisely, the question of interest is to analyze the effect of
letting $\varepsilon_1$ and $\varepsilon_2$ go independently to 0 in
the expression $(f\star \varphi_{\varepsilon_1,\varepsilon_2})(\bx)$.
We prove in particular (Theorem \ref{proannex2}) that $(f\star \varphi_{\varepsilon_1,\varepsilon_2})(\bx)\to f(\bx)$ for $\lambda_n$-almost all $\bx$ if $f$ and $\varphi$ are suitably restricted. 
\medskip

The issues discussed in the present appendix fall within the field of maximal functions and approximation of the identity \citep[][]{stein,Wheeden}. The novelty is that we allow the family $\{\varphi_{\varepsilon_1,\varepsilon_2}:\varepsilon_1>0,\varepsilon_2>0\}$ (the so-called approximation of the identity) to depend upon {\bf two independent} parameters $\varepsilon_1$ and $\varepsilon_2$. Interestingly, the real analysis literature offers little help with respect to this important question, which is however fundamental in the study of multivariate nonparametric estimates. Valuable ideas and comments in this respect are included in \citet{DK}.
\medskip

Let $\varphi$ be an integrable function on $\mathbb R^n=\mathbb R^{n_1}\times\mathbb R^{n_2}$, termed ``the kernel'' hereafter. It is assumed throughout that $\varphi$ is a product kernel, of the form 
\begin{equation}
\label{kerprod}
\varphi(\bx) = \varphi_1(\bx_1)\varphi_2(\bx_2), \quad \bx=(\bx_1,
\bx_2) \in \mathbb R^{n_1}\times \mathbb R^{n_2}.
\end{equation} 
For $\varepsilon_1>0$ and $\varepsilon_2>0$, we set
$$
\varphi_{\varepsilon_1,\varepsilon_2}(\bx)=\frac{1}{\varepsilon_1^{n_1}\varepsilon_2^{n_2}}
\varphi_1\left (\frac{\bx_1}{\varepsilon_1}\right)\varphi_2\left (\frac{\bx_2}{\varepsilon_2}\right).
$$
We will need the following assumption: 
\medskip

\noindent\textbf{Assumption [K]}\quad For $i=1,2$, the functions
$$\psi_{i}(\bx_i) = \sup_{\|\mathbf y_i\|\geq \|\bx_i\|} \left |\varphi_i(\mathbf y_i)\right|,\quad \bx_i \in \mathbb R^{n_i},$$ are in $L^1(\mathbb R^{n_i})$, with 
$$\int_{\mathbb R^{n_i}} \psi_i(\bx_i)\mbox{d}\bx_i \leq \sqrt{A}<\infty.$$

If $f$ is a locally integrable function in $\mathbb R^n$, we also denote by $M_{12}f$ the associated Hardy-Littlewood maximal function with two
degrees of freedom. It is defined for $\bx=(\bx_1,\bx_2)$ by
$$(M_{12}f)(\bx) = \sup_{\varepsilon_1, \varepsilon_2 >0}
\left [\frac{1}{V_{n_1}\varepsilon_1^{n_1}V_{n_2}\varepsilon_2^{n_2}}
\int_{\mathcal B_{n_1}(\bx_1,\varepsilon_1)} \int_{\mathcal B_{n_2}(\bx_2,\varepsilon_2)}
\left |f(\mathbf y_1,\mathbf y_2)\right|\mbox{d}\mathbf y_1\mbox{d}\mathbf y_2\right],$$
where $\mathcal B_{n_1}(\bx_1,\varepsilon_1)$ (respectively, $\mathcal
B_{n_2}(\bx_2,\varepsilon_2)$) is the closed ball in $\mathbb  R^{n_1}$
(respectively, $\mathbb R^{n_2})$, with center at $\bx_1$
(respectively, $\bx_2$) and radius $\varepsilon_1$ (respectively,
$\varepsilon_2$), and $V_{n_1}$ (respectively, $V_{n_2}$) is the
volume of the unit ball in $\mathbb R^{n_1}$ (respectively, $\mathbb R^{n_2}$).
\medskip

Our objective is to prove the following theorem, which
is a more general version of Theorem \ref{WZ2}.
\begin{theo}
\label{proannex2}
Let $f$ be a measurable function in $\mathbb R^n$ satisfying 
\begin{equation}  \label{mcnab} \int_{\mathbb R^n} \left
    |f(\bx)\right|\left(1+\log^+\left
      |f(\bx)\right|\right)\emph{d}\bx<\infty,\end{equation} 
and let $\varphi$ be a product kernel of the form (\ref{kerprod}) satisfying Assumption $[{\bf K}]$. Assume, in addition, that $\int_{\mathbb R^n} \varphi(\bx) \emph{d}\bx =1$.
\begin{itemize} 
\item[$(i)$] For $\lambda_n$-almost all $\bx\in\mathbb R^n$, $\lim_{\varepsilon_1,\varepsilon_2\to0}
(f\star\varphi_{\varepsilon_1,\varepsilon_2})(\bx) = f(\bx)$.
\item[$(ii)$] For $\lambda_n$-almost all $\bx\in\mathbb R^n$, 
$$
\sup_{\varepsilon_1,\varepsilon_2 >0} \left|(f \star
\varphi_{\varepsilon_1,\varepsilon_2})(\bx)\right| \leq A(M_{12}f)(\bx)<\infty,
$$
where $A$ is the constant of Assumption $[{\bf K}]$.
\item[$(iii)$] Moreover, if $f$ is in $L^q(\mathbb R^n)$, $1<q\leq \infty$, then 
$M_{12} f$ is  in $L^q(\mathbb R^n)$
and 
$$\|M_{12}f\|_q\leq c_q\|f\|_q,$$
where the constant $c_q$ depends only on $q$ and the dimension $n$.
\end{itemize} 
\end{theo}

\noindent{\bf Proof of Theorem \ref{proannex2}}\quad 
To prove the theorem, we will need some general results on singular
integrals and Hardy-Littlewood maximal functions.  As shown in page 50 of \citet{guzman}, for all $\alpha>0$ and a locally
integrable $f$,
\begin{equation}
\label{eq.guzman}
\lambda_n\left (\{ \bx\in\mathbb R^{n}: (M_{12}f)(\bx)>\alpha\}\right)
\leq
c\int_{\mathbb R^{n}} \frac{\left |f(\bx)\right|}{\alpha}\left (1+\log^+\frac{\left |f(\bx)\right|}{\alpha}\right) \mbox{d}\bx,
\end{equation}
where $c$ is a constant independent of $f$ and $\alpha$. 
This result will be crucial in our proof. 
It easily follows that whenever
$$\int_{\mathbb R^{n}}
\left |f(\bx)\right|\left(1+\log^+{\left |f(\bx)\right|}\right) \mbox{d}\bx <\infty,$$
then $(M_{12}f)(\bx)<\infty$ at $\lambda_n$-almost all $\bx$.
\paragraph{Proof of $(ii)$} The proof follows arguments of pages 63-64 of \citet{stein}.  For $i=1,2$, with a slight abuse of notation, we write
$\psi_i(r_i)=\psi_i(\bx_i)$ if $r_i=\|\bx_i\|$. This should cause no confusion since each $\psi_i$ is anyway radial. Observe that, for $i=1,2$,
$$
\int_{r_i/2\leq \|\bx_i\|\leq r_i} \psi_i(\bx_i) \mbox{d}\bx_i \geq \psi_i(r_i)
\int_{r_i/2\leq \|\bx_i\|\leq r_i} \mbox{d}\bx_i \propto \psi_i(r_i) r^{n_i}_i.
$$
Therefore, the assumption $\psi_i \in L^1(\mathbb R^{n_i})$ proves that $r_i^{n_i}
\psi_i(r_i) \to 0$, as $r_i\to 0$ or  $r_i\rightarrow \infty$. To prove $(ii)$, it is enough to show that for all nonnegative $f$ satisfying (\ref{mcnab}), all $\varepsilon_1>0, \varepsilon_2>0$,
\begin{equation}
\label{16}
(f\star\psi_{\varepsilon_1,\varepsilon_2})(\bx) \leq A (M_{12}f)(\bx),
\end{equation}
where 
$$\psi_{\varepsilon_1,\varepsilon_2}(\bx)=\frac{1}{\varepsilon_1^{n_1}\varepsilon_2^{n_2}}
\psi_1\left (\frac{\bx_1}{\varepsilon_1}\right)\psi_2\left(\frac{\bx_2}{\varepsilon_2}\right), \quad \bx=(\bx_1,\bx_2)\in \mathbb R^n.$$
Set $\psi=\psi_1\psi_2$. Since assertion (\ref{16}) is clearly translation invariant (with respect to $f$) and also dilatation invariant (with respect to $\psi$), it suffices to show that 
$$
(f \star\psi)(\mathbf 0) \leq A(M_{12}f)(\mathbf 0).
$$
Moreover, recalling (\ref{eq.guzman}), we may clearly assume that $(M_{12}f)(\mathbf 0)<\infty$. For $i=1,2$, denote by $S^{n_i-1}$ the unit $(n_i-1)$-sphere in $\mathbb R^{n_i}$ and let $\sigma_i$ be the corresponding spherical measure. We set as well
\begin{align*}
\ell(r_1,r_2)&=\int_{S^{n_1-1}} \int_{S^{n_2-1}} f(r_1 \bx_1, r_2 \bx_2)
\mbox{d}\sigma_1(\bx_1) \mbox{d}\sigma_2(\bx_2),\\
\Lambda_1(r_1,r_2) &= \int_0^{r_1} \ell(u_1,r_2) u_1^{n_1-1}
\mbox{d}u_1= \int_0^{r_2} \Lambda_1(r_1,u_2) u_2^{n_2-1} \mbox{d}u_2 \\
&=\int_0^{r_1}\int_0^{r_2} \ell(u_1,u_2) u_1^{n_1-1} u_2^{n_2-1}
\mbox{d}u_1 \mbox{d}u_2,
\end{align*} 
and will repeatedly use the inequality
\begin{equation}
\label{1206}
\Lambda(r_1,r_2)=\int_{\mathcal B_{n_1}(\mathbf 0,r_1)} \int_{\mathcal B_{n_2}(\mathbf 0,r_2)}f(\bx)\mbox{d}\bx\leq V_{n_1}.V_{n_2}r_1^{n_1} r_2^{n_2} (M_{12}f)(\mathbf 0).
\end{equation}
\medskip

With this notation, we have
\begin{align*}
(f\star\psi)(\mathbf 0)&= \int_0^\infty \int_0^\infty \ell(r_1,r_2)\psi_1(r_1)r_1^{n_1-1}\psi_2(r_2)r_2^{n_2-1} \mbox{d}r_1 \mbox{d}r_2 \\
&=\lim_{{\tiny \begin{array}{c} \varepsilon_1\to 0\\
N_1\to \infty\\\varepsilon_2\rightarrow 0\\N_2\to
\infty \end{array} }}
\int_{\varepsilon_2}^{N_2} \left[\int_{\varepsilon_1}^{N_1}
  \ell(r_1,r_2) \psi_1(r_1) r_1^{n_1-1} \mbox{d}r_1\right] \psi_2(r_2)
r_2^{n_2-1} \mbox{d}r_2.
\end{align*} 
Denote by $I_1(\varepsilon_1,N_1)$ the integral inside the
brackets. We may write, using an integration by parts (in the sense of Stieltj\`es-Lebesgue),
$$
I_1(\varepsilon_1,N_1) = \int_{\varepsilon_1}^{N_1} \Lambda_1(r_1,r_2)
\mbox{d}\left(-\psi_1(r_1)\right) + \Lambda_1(N_1,r_2)\psi_1(N_1) -
\Lambda_1(\varepsilon_1, r_2)\psi_1(\varepsilon_1).
$$
Consequently,
\begin{align*}
&\int_{\varepsilon_2}^{N_2} I_1(\varepsilon_1,N_1)\psi_2(r_2)
r_2^{n_2-1} \mbox{d}r_2 = I_A + I_B - I_C\\
&\quad = \int_{\varepsilon_2}^{N_2}\int_{\varepsilon_1}^{N_1}
\Lambda_1(r_1,r_2) \mbox{d}\left (-\psi_1(r_1)\right) \psi_2(r_2)r_2^{n_2-1} \mbox{d}r_2\\
&\qquad +
\int_{\varepsilon_2}^{N_2} \Lambda_1(N_1,r_2)\psi_1(N_1)\psi_2(r_2)
r_2^{n_2-1} \mbox{d}r_2\\
&\qquad -
\int_{\varepsilon_2}^{N_2}
\Lambda_1(\varepsilon_1,r_2)\psi_1(\varepsilon_1) \psi_2(r_2)
r_2^{n_2-1} \mbox{d}r_2.
\end{align*} 
Each term of the sum is analyzed separately. Using again an integration by parts, we are led to
\begin{align*}
I_A&=\int_{\varepsilon_1}^{N_1}\Big[ \int_{\varepsilon_2}^{N_2}
  \Lambda(r_1,r_2) \mbox{d}\left (-\psi_2(r_2)\right)+\Lambda(r_1,N_2)\psi_2(N_2)\\
  &\qquad \qquad -\Lambda(r_1,\varepsilon_2)\psi_2(\varepsilon_2)\Big]\mbox{d}\left(-\psi_1(r_1)\right)\\
&=
\int_{\varepsilon_1}^{N_1}\int_{\varepsilon_2}^{N_2} \Lambda(r_1,r_2)
\mbox{d}\left (-\psi_1(r_1)\right)\mbox{d}\left (-\psi_2(r_2)\right)\\
&\quad +
\int_{\varepsilon_1}^{N_1} \Lambda(r_1,N_2)\psi_2(N_2)\mbox{d}\left (-\psi_1(r_1)\right)
\\
&\quad - 
\int_{\varepsilon_1}^{N_1}
\Lambda(r_1,\varepsilon_2)\psi_2(\varepsilon_2)\mbox{d}\left (-\psi_1(r_1)\right) \\
&=
A_1 + A_2 - A_3.
\end{align*} 
The main term, $A_1$, is handled as follows via inequality (\ref{1206}):
\begin{align*}
A_1&\leq V_{n_1}.V_{n_2}(M_{12}f)(\mathbf 0) \int_0^\infty \int_0^\infty r_1^{n_1}
r_2^{n_2 } \mbox{d}\left(-\psi_1(r_1)\right)\mbox{d}\left(-\psi_2(r_2)\right)\\
& \leq A(M_{12}f)(\mathbf 0)
\end{align*}
since for $i=1,2$, we have
$$V_{n_i}\int_0^\infty r_i^{n_i}\mbox{d}\left(-\psi_i(r_i)\right)=\int_{\mathbb R^{n_i}} \psi_i(\bx_i)\mbox{d}\bx_i\leq\sqrt{A},$$
by Assumption $[\mathbf{K}$]. The remaining terms, $A_2$ and $A_3$, converge to $0$. To see this, just note that
$$
A_2  \leq V_{n_1}.V_{n_2}(M_{12}f)(\mathbf 0) \times N_2^{n_2}\psi_2(N_2)\int_{0}^{\infty}
r_1^{n_1} \mbox{d}\left (-\psi_1(r_1)\right),
$$
which goes to $0$ since the integral is convergent and
$N_2^{n_2}\psi_2(N_2) \to 0$ as $N_2\to \infty$.  Similarly,
$$A_3\leq V_{n_1}.V_{n_2} (M_{12} f)(\mathbf 0)\times \varepsilon_2^{n_2}\psi_2(\varepsilon_2)
\int_0^\infty r_1^{n_1}\mbox{d}\left (-\psi_1(r_1)\right).
$$
The term on the right-hand side tends to 0 since $\varepsilon_2^{n_2}\psi_2(\varepsilon_2)\to 0$ as $\varepsilon_2\to 0$. Using similar arguments, it is easy to prove that $I_B$ and $I_C$ go to 0 as $\varepsilon_1, \varepsilon_2\to 0$ and $N_1,N_2\to \infty$. Proof of $(ii)$ is therefore complete.
\paragraph{Proof of $(i)$} For the sake of clarity, the proof is divided into three steps.
\medskip

\noindent{\bf Step 1}\quad If $f$ is continuous and has compact support, then the result is easy to verify. Indeed, we have in this case
$$
(f\star\varphi_{\varepsilon_1,\varepsilon_2})(\bx) = \int_{\mathbb R^{n_1}}\int_{\mathbb  R^{n_2}}
f(\bx_1-\varepsilon_1 \mathbf y_1, \bx_2-\varepsilon_2 \mathbf y_2) \varphi(\mathbf y_1,\mathbf y_2) \mbox{d}\mathbf y_1
\mbox{d}\mathbf y_2,
$$
whence, using the fact that $\int_{\mathbb R^n} \varphi(\bx) \mbox{d}\bx =1,$
\begin{align*}
&\left |(f\star\varphi_{\varepsilon_1,\varepsilon_2})(\bx) -f(\bx)\right|\\
&\quad \leq \int_{\mathbb R^{n_1}}\int_{\mathbb R^{n_2}}
\left |f(\bx_1-\varepsilon_1 \mathbf y_1, \bx_2-\varepsilon_2 \mathbf y_2) -f(\bx)\right|.\left |\varphi(\mathbf y_1,\mathbf y_2)\right| \mbox{d}\mathbf y_1
\mbox{d}\mathbf y_2 \\
&\quad \leq \sup_{\bx_1,\bx_2,\mathbf y_1,\mathbf y_2} \left |f(\bx_1-\varepsilon_1 \mathbf y_1,
\bx_2-\varepsilon_2 \mathbf y_2) -f(\bx)\right| \int_{\mathbb R^{n_1}}\int_{\mathbb R^{n_2}}
 \left |\varphi(\mathbf y_1,\mathbf y_2)\right| \mbox{d}\mathbf y_1
\mbox{d}\mathbf y_2.
\end{align*}
Since $f$ is uniformly continuous, this term tends to $0$. 
\medskip

\noindent{\bf Step 2}\quad We establish that $\lim_{\varepsilon_1,\varepsilon_2\to 0} (f\star \varphi_{\varepsilon_1,\varepsilon_2})(\bx)$ exists for $\lambda_n$-almost all $\bx \in \mathbb R^n$. As for now, to ease the notation, we set $g^{\star}_{\varepsilon_1,\varepsilon_2}(\bx)=(g \star \varphi_{\varepsilon_1,\varepsilon_2})(\bx)$, and let
$$
(\Omega g)(\bx) = \left |\limsup_{\varepsilon_1,\varepsilon_2\to 0}
g^{\star}_{\varepsilon_1,\varepsilon_2}(\bx) -
\liminf_{\varepsilon_1,\varepsilon_2\rightarrow 0} 
g^{\star}_{\varepsilon_1,\varepsilon_2}(\bx)\right|.
$$
Let $\alpha>0$ and $\delta>0$ be arbitrary. Thanks to Proposition~\ref{lemannex1} at the end of the section, we may write
$f=h+g$, where $h$ is continuous with compact support and $g$ is such that
$$
\int_{\mathbb R^n}\frac{\left |g(\bx)\right|}{\alpha}\left (1+\log^+ \frac{\left |g(\bx)\right|}{\alpha}\right) \mbox{d}\bx \leq \delta.
$$
By $(ii)$, we have at $\lambda_n$-almost all $\bx$, $(\Omega g)(\bx) \leq 2A(M_{12} g)(\bx)$. Thus, by (\ref{eq.guzman}),
$$\lambda\left (\{ \bx \in \mathbb R^n: (\Omega g)(\bx)>2A\alpha\}\right)\leq
c\int_{\mathbb R^n}\frac{\left |g(\bx)\right|}{\alpha}\left (1+\log^+ \frac{\left |g(\bx)\right|}{\alpha}\right)\mbox{d}\bx  \leq c\delta.$$ 
Clearly, $\Omega f\leq\Omega g + \Omega h$ and, by Step 1, $\Omega h
\equiv 0$. Therefore  
$$
\lambda\left (\{\bx \in \mathbb R^n: (\Omega f)(\bx)>2A\alpha\}\right)\leq c\delta.
$$
Since $\alpha$ and $\delta$ are arbitrary, we conclude that $\lambda\left (\{\bx \in \mathbb R^n: (\Omega f)(\bx)>0\}\right)=0$.
\medskip
 
\noindent{\bf Step 3}\quad  We finally prove that, for $\lambda_n$-almost all $\bx \in \mathbb R^n$,
$$\lim_{\varepsilon_1,\varepsilon_2\to 0}
f^{\star}_{\varepsilon_1,\varepsilon_2}(\bx)=f(\bx).$$
Set $f_1(\bx) = \lim_{\varepsilon_1,\varepsilon_2\rightarrow
  0} f^{\star}_{\varepsilon_1,\varepsilon_2}(\bx)$ (this limit exists $\lambda_n$-almost everywhere by Step 2). Fix $\alpha>0$, $\delta>0$, and choose $h$ continuous with compact support as in Step 2 such that
  $$
\int_{\mathbb R^n}\frac{\left |(f-h)(\bx)\right|}{\alpha}\left (1+\log^+ \frac{\left |(f-h)(\bx)\right|}{\alpha}\right) \mbox{d}\bx \leq \delta.
$$
For $\lambda_n$-almost all $\bx\in\mathbb R^n$,
$$\left |f(\bx)-f_1(\bx)\right|\leq\left |f(\bx)-h(\bx)\right| +
|\lim_{\varepsilon_1,\varepsilon_2\rightarrow 0}
h^{\star}_{\varepsilon_1,\varepsilon_2}(\bx) -
\lim_{\varepsilon_1,\varepsilon_2\rightarrow 0}
f^{\star}_{\varepsilon_1,\varepsilon_2}(\bx) |= A_1 + A_2.$$
By $(ii)$, 
$$A_2\leq \sup_{\varepsilon_1,\varepsilon_2>0}
\left |(f-h)^{\star}_{\varepsilon_1,\varepsilon_2}(\bx)\right| \leq A\left (M_{12}|f-h|\right)(\bx).$$
Thus,
\begin{align*}
&\lambda\left (\{\bx\in \mathbb R^n: \left |f(\bx)-f_1(\bx)\right|>2A \alpha\}\right) \\
&\quad \leq \lambda\left (\{\bx\in \mathbb R^n: \left |f(\bx)-h(\bx)\right|>A \alpha\}\right) \\
& \qquad + \lambda\left (\{\bx\in \mathbb R^n: \left (M_{12}|f-h|\right)(\bx)>\alpha\}\right)
\\
&\quad \leq \frac{\|f-h\|_1}{A\alpha} + c \int_{\mathbb R^n} \frac{\left |(f-h)(\bx)\right|}{\alpha}\left (1+\log^+ \frac{\left |(f-h)(\bx)\right|}{\alpha}\right) \mbox{d}\bx \\
&\quad \leq\left(\frac{1}{A}+c\right)\delta.
\end{align*} 
In the second inequality, we used Markov's inequality together with inequality (\ref{eq.guzman}). Since both $\alpha$ and $\delta$ can be chosen arbitrarily, we conclude
that 
$$\lambda\left (\{\bx \in \mathbb R^n: \left |f(\bx)-f_1(\bx)\right|>0\}\right)=0.$$
\paragraph{\bf Proof of $(iii)$} The proof is adapted from page 307 of \citet{Zygmund}. Let the partial maximal functions be defined for $\bx=(\bx_1,\bx_2)$ by
$$
(M_1 f)(\bx) = \sup_{\varepsilon_1>0}
\left[\frac{1}{V_{n_1}\varepsilon_1^{n_1}} \int_{\mathcal B_{n_1}(\bx_1,\varepsilon_1)}
|f(\mathbf y_1,\bx_2)| \mbox{d}\mathbf y_1\right]
$$
and
$$
(M_2 f)(\bx) = \sup_{\varepsilon_2>0}\left [
\frac{1}{V_{n_2}\varepsilon_2^{n_2}} \int_{\mathcal B_{n_2}(\bx_2,\varepsilon_2)}
|f(\bx_1,\mathbf y_2)| \mbox{d}\mathbf y_2\right].
$$
From these definitions, it is clear that $(M_{12} f)(\bx)\leq \left (M_1(M_2 f)\right)(\bx)$. But, for $1<q\leq \infty$, $f_1\in L^q(\mathbb R^{n_1})$, $f_2\in L^q(\mathbb R^{n_2})$, it is known \citep[see, e.g.,][Theorem 1, page 5]
{stein}, that
$$
\left \|M_1f\right\|_q \leq c_{1,q} \left \|f_1\right\|_q \mbox{ and } \left \|M_2 f\right\|_q \leq c_{2,q} \|f_2\|_q,
$$
where the constants $c_{1,q}$ and $c_{2,q}$ depend only on $n_1$, $n_2$ and $q$. It immediately follows that  $\left\|M_{12}f\right\|_q^q  \leq  c_{1,q}^q c_{2,q}^q \|f\|_q^q$. This concludes the proof of the theorem. 
\hfill $\blacksquare$
\medskip

\begin{pro}
\label{lemannex1}
Let $\Phi : \mathbb R^+ \rightarrow \mathbb R^+$ be a continuous and nondecreasing function satisfying $\Phi(0)=0$, and let $f$ be a measurable function from $\mathbb R^n$ to $\mathbb R$ such that $\int_{\mathbb R^n} \Phi\left (\left |f(\bx)\right|\right) \emph{d}\bx < \infty$. Then, for all $\delta>0$, there exists a function $h$ continuous with compact support
such that 
$$\int_{\mathbb R^n} \Phi\left (\left |f(\bx)-h(\bx)\right|\right) \emph{d}\bx \leq \delta.$$
\end{pro}
\noindent{\bf Proof of Proposition \ref{lemannex1}}\quad First, assume that $f(\bx)\geq 0$ for all $\bx$. Take $\{f_t\}$ a sequence of nonnegative continuous functions, each with compact support and such that $0\leq f_t(\bx) \uparrow f(\bx)$ at $\lambda_n$-almost all $\bx \in \mathbb R^n$. For such an $\bx$, by the continuity of $\Phi$ at 0, one has $\Phi(f(\bx)-f_t(\bx))\to \Phi(0)=0$. Since $\Phi(f(\bx)-f_t(\bx))\leq \Phi(f(\bx))$ and $\Phi(f)$ is in $L^1(\mathbb R^n)$ by assumption,  we may apply Lebesgue's dominated convergence theorem and conclude that
$$\int_{\mathbb R^n} \Phi\left (f(\bx)-f_t(\bx)\right) \mbox{d}\bx \to 0 \quad \mbox{as }t \to \infty.$$If we drop the assumption that $f(\bx)\geq 0$, we may split $f$ into positive and negative part and apply the above result.
\hfill $\blacksquare$
\paragraph{Acknowledgments} We thank two anonymous referees for valuable comments and insightful suggestions. 

\bibliography{biblio-abc}
\end{document}